\documentclass{amsart}
\usepackage{amssymb}
\usepackage{amsfonts}
\usepackage{amsmath}
\usepackage[legalpaper,bookmarks=true,colorlinks=true,linkcolor=blue,citecolor=blue]{hyperref}
\usepackage{graphicx}%
\usepackage{fancyhdr}
\usepackage{color}
\usepackage[mathlines]{lineno}
\usepackage{lscape}
\usepackage{epsfig}
\usepackage{natbib}
\usepackage{geometry}
\usepackage{tgbonum}
\fontfamily{qcr}\selectfont

\newtheorem{theorem}{Theorem}
\theoremstyle{plain}

\newtheorem{lemma}{Lemma}

\newtheorem{proposition}{Proposition}
\newtheorem{remark}{Remark}

\numberwithin{equation}{section}

\newcommand{\Bin}{\bigskip \noindent}

\newcommand{\Ni}{\noindent}

\begin{document}
\Large
\title[Moving Central limit theorems]{Central limit theorems for associated possibly moving partial sums and application to the non-stationary invariance principle}
\author{Akim Adekpedjou}
\author{Aladji Babacar Niang}
\author{Ch\'erif Mamadou Moctar Traor\'e}
\author{Gane samb Lo}

\begin{abstract} General Central limit theorem deals with weak limits (in type) of sums of row-elements of array random variables. In some situations as in the invariance principle problem, the sums may include only parts of the row-elements. For strictly stationary arrays (stationary for each row), there is no change to the asymptotic results. But for non-stationary data, especially for dependent data, asymptotic laws of partial sums moving in rows may require extra-conditions to exist. This paper deals with central limit theorems with Gaussian limits for non-stationary data. Our main focus is on dependent data, particularly on associated data. But the non-stationary independent data is also studied as a learning process. The results are applied to finite-distributional invariance principles for the types of data described above. In Moreover, results for associated sequences are interesting and innovative. Beyond their own interest, the results are expected to be applied for random sums of random variables and next in statistical modeling in many disciplines, in Actuarial sciences for example. \\

\noindent Aladji Babacar Niang\\
LERSTAD, Gaston Berger University, Saint-Louis, SENEGAL\\
Email: niang.aladji-babacar@ugb.edu.sn, aladjibacar93@gmail.com\\

\noindent Cherif Mamadou Moctar TRAORE\\
LMA, FST, EDSTM, University of  Technical and Technological Sciences of Bamako (USTTB), MALI.\\
LERSTAD, Gaston Berger University (UGB), Saint-Louis, SENEGAL.\\
Email :traore.cherif-mamadou-moctar@ugb.edu.sn, cheriftraore75@yahoo.com\\

\noindent $^{\dag}$ Gane Samb Lo.\\
LERSTAD, Gaston Berger University, Saint-Louis, SENEGAL (main affiliation).\newline
LSTA, Pierre and Marie Curie University, Paris VI, FRANCE.\newline
AUST - African University of Science and Technology, Abuja, NIGERIA\\
gane-samb.lo@edu.ugb.sn, gslo@aust.edu.ng, ganesamblo@ganesamblo.net\\
Permanent address : 1178 Evanston Dr NW T3P 0J9,Calgary, Alberta, CANADA.\\


\end{abstract}
\maketitle

\section{Introduction} \label{sec_01}

\Ni Moving partial sums are closely related to invariance principles, which in turn play an important role in many areas of applications such as Finance, Actuarial Sciences, Demography, etc. As an example, consider the claims problem for an insurer, whose clients subscribe to specific products through determined policies. In vehicule insurance for instance, the policy may include that at each accident, the client makes a claim $X$, which depends on many factors as the severity of the crash for instance. For simplicity, we suppose that the claims are reported at discrete times $nt_0$, where $t_0$ is a fixed period of time that may in days, weeks or months. At each time  $jt_0$, the claim is a random variables $X_j$. So, at time $nt_0$, the total claim (referred as the total loss) up to time $nt_0$ is given by the equation

\begin{equation}
S_n=\sum_{1\leq j\leq n} X_{j}. \label{totLossD}
\end{equation}

\Bin The insurer should have a accurate estimation of $S_n$ to fix the premiums by clients should pay at the establishment of the policies, otherwise the ruin would be highly probable. We remark that the the discrete time modeling of claims \eqref{totLossD}  can be extended to a continuous time one. In such a case, the number of reported claims up to time $t$, say $N(t)$, is a random variable and the total loss up to $t$ is 

$$
S_{t}=\sum_{j=1}^{N(t)} X_j.
$$ 

\Bin Now, we suppose that the insurer has a capital $u$ at the beginning, and that the premiums can be linearized, say as $ct$, the surplus process (measuring the financial balance of the insurer) at time $t$ can be given as

$$
P_{t}=u + ct - S_{t}. 
$$

\Bin Although though the model uses continuous time, in practice, the time is discretized into multiples of a unit of time $t_0>0$ and the model becomes, for $t_n=nt_0$, $n\geq 0$,

$$
P_{t_n}=u + c t_n - \sum_{j=1}^{N(t_n)} X_j=: u + ct_n - c_n \frac{\sum_{j=1}^{N(t_n)} X_j}{c_n}, \ n\geq 0. 
$$

\Bin Finding the limiting law of the stochastic process 

$$
\left\{Y_n(t)=\frac{\sum_{j=1}^{N(t_n)} X_j}{c_n}, \ 0\leq t\leq T, \ n\geq 1 \right\},
$$

\Bin for $T>0$, for an appropriate sequence of normalization coefficients $(c_n)_{n\geq 1}$, to a stochastic process $\{Y(t), \ 0\leq t\leq T\}$ is the essence of the invariance principle problem (or functional central limit theorem) problem. For independent data, the most used limiting law is a Brownian motion. But, even in that case,
 the general solution is a L\'evy process $Y$.\\ 

\Ni Usually, $N(t)$ is taken as a Poisson process. It is reasonable to expect that at least, the insurer should avoid incurring a ruin, say at a time $t_{n(r)}$, such that $P_{t_{n(r)}}<0$. An approximation of the probability ruin is given by

\begin{equation}
P_{t_{n(r)}} \approx u + ct_{n(r)} - Y(t_{n(r)}), \ n\geq 0, \label{Wealth} 
\end{equation}

\Bin and

$$
n(r)= \inf \{n\geq 1,  P_{t_n}<0\}.
$$

\Bin If $Y$ is accurately estimated, Equation \eqref{Wealth} may help in pre-setting $c$ and $u$ before contracting policies, to ensure profit and avoid ruin. To learn more on such modeling, the reader is directed to \cite{klugman} page 252, \cite{grandel91} and the references therein. For independent and square integrable data, the class of possible weak limits is exactly that of of infinitely decomposable laws and the associated invariance principle leads to Levy processes, the Brownian motion and the Poisson process being among them (see \cite{loeve}, \cite{applebaum}, \cite{niang-summands-01}, etc.).\\

\Ni The general problem of finding the weak law, in the random scheme as

$$
S_n(t)=\sum_{j=1}^{N([nt]} X_j, \ 0\leq t\leq T, \ n\geq 1
$$

\Bin and, in the non-random scheme as

$$
S_n(t)=\sum_{j=1}^{[nt]} X_j, \ 0\leq t\leq T, \ n\geq 1,
$$

\Bin (usually for $T=1$) is the core of the invariance principle problem. This problem is hard and quite general since we do not necessarily know the dependence type between the losses $X_j$'s nor do we always have that $N(\circ)$ is a Poisson process. The  main results and achievements in the literature are obtained for independent data and when $N$ is a classical Poisson process.\\  

\Bin The problem of moving partial sums arise in the important setting of functional weak limits. In what follows, we provide some background. Consider a sequence of centered random variables $\left(X_{n}\right)_{n\geq 1}$ defined on the same probability space $(\Omega,\mathcal{A},\mathbb{P})$. Let $\left(p(n)\right)_{n\geq 1}$ be an arbitrary sequence of positive integers. We define

\begin{equation}
S^{\prime}_0=0,  \ S_{n}^{\prime}=\sum_{k=p(n)+1}^{n+p(n)} X_k \,\,\,\,\, and \,\,\,\,\, s_{n}^{\prime2}=\mathbb{V}ar \left(S_{n}^{\prime}\right), n\geq 1,\label{mpsum-not}
\end{equation}

\bigskip \noindent with, for all $k\geq 1$,  $\sigma_{k}^{2} = \mathbb{E}(X_{k}^{2})$ and $F_k(x)=\mathbb{P}(X_k\leq x)$, $x \in \mathbb{R}$.\\

\noindent If $p(n)=0$ for all $n\geq 1$, we find ourselves in studying the usual partial sums $S_{n} = X_{1} + X_{2} + ... + X_{n}$ and the partial sums of variances $s_{n}^{2} = \mathbb{V}ar(S_n)$.\\

\noindent In some situations, we may be concerned not with all the partial sums $S_n$ beginning by the first \textit{r.v.} $X_1$ but with partial sums that can start at any part of the sequence $\left(X_{n}\right)_{n\geq 1}$. In \eqref{mpsum-not}, the partial sums $S^{\prime}_n$ begin with the \textit{r.v.} $X_{p(n)+1}$ and is the sum of all the $n$ observations with index greater that $p(n)$. A more appropriate notation should be $S_{n,p(n)}$ that we denote as $S^{\prime}_n$, given that the sequence $(p(n))_{n\geq 1}$ is already defined. For example, when dealing with invariance principles, we need to have the limit of the sequence of stochastic processes 
$\left\{Y_n(t), \ 0\leq t\leq T, \ n\geq 1\right\}$ to some stochastic process $\left\{Y(t), \ 0\leq t\leq T\right\}$. The state of the art (see \cite{billingsley}, \cite{vaart}, \cite{ips-wcia-ang}, etc.) expresses that, in order that weak convergence holds, we need to have the convergence in finite-distribution and that the sequence is uniformly tight (as in \cite{billingsley}) or asymptotically tight (as in \cite{vaart}). For now, we focus on the weak convergence in finite distributions, i.e., of the vectors 

$$
(Y_n(t_j))_{\ 1\leq j \leq k}=\biggr(\frac{S_{[nt_j]}-S_{[nt_{j-1}]}}{s_n}\biggr)_{1\leq j \leq k},
$$

\Bin for $0=t_0<t_1<\cdots<t_k=T$ (usually $T=1$), $k\geq 2$. It is clear that we have, for each $j$-th component:

$$
Y_n(t_j)=\frac{S^{\prime}_{[nt_j]-[nt_{j-1}]}}{s_n}, \ with \ p(n)=[nt_{j-1}], \ 1\leq j \leq k.
$$

\Bin Here the relation between $s_n$ and $s_n^{\prime}$ plays a major role as we will see later. We refer to such partial sums as moving partial sums, since the sequence $\left(p(n)\right)_{n\geq 1}$ is arbitrary. As we will see, the handling weak laws of moving partial sums (\textit{MPS}) can be a lot easier for independent \textit{r.v.}'s. Nevertheless, we have to make sure that all steps are rigorously taken into account. \textit{But, for dependent data, the situation requires more attention and can get very complicated. The problem is even more serious if the sequence $(p(n))_{n\geq 1}$ is random as expected in applications, especially in Actuarial Sciences and Finance}.\\

\Ni In this paper, we focus on associated data introduced in \ref{sec_03_subsec_01} as associated data. We will assume knowledge of basic definitions and results on assiciated data, thus, we refer  the reader to \cite{GH-gslo} for a quick review. In \ref{sec_03_subsec_01}, we will review the most important facts on associated data. For much more details, \cite{rao} and \cite{bulinski} are more appropriate. Nevertheless, we will present a thorough review of some important results on independent data, in order to facilitate the passage to dependent data.\\

\Ni Hence, the purpose of this paper is to extend Gaussian central limits for independent data to Gaussian limits for \textit{MPS} and to compare the classical conditions (Lyapounov and Lynderberg, Uniform Asymptotic Negligibility (\textit{UAN}), Bounded Variance Hypothesis (\textit{BVH}), Convergence Variance Hypothesis (\textit{CVH}), etc.) for full sums and moving partial sums. We will then study how such results can be used in invariance principles to re-scaled Brownian motions. Since handling invariance principles for dependent data, here for associated data, requires similar asymptotic weak laws for \textit{MPS} and their applications, this will lead to more general results in invariance principles for associated data in comparison with current achievements in particular in \cite{GH-paulo}. The present analysis will open the door to studying asymptotic weak laws for \textit{MPS} for other types of dependence and for other specific types limit laws, i.e., for any infinitely decomposable type limiting laws, in particular Poisson laws. Future research works will focus on generalizations of the results to random sums and their data-driven applications.\\

\Ni  The rest of the paper is organized as follows. In Section \ref{sec_02}, we study the \textit{MPS} in the independence situation under the condition 

$$
\forall t \in (0,1), \ s^2_{[nt]}/s^2_n \rightarrow a(t) \ as \ n\rightarrow +\infty, 
$$

\Bin where $a(t)$ is a non-decreasing function in $t\in [0,1]$. Next, we apply the results to finite-distributional invariance principles with weak convergence to re-scaled Brownian motions. We get moving versions for Lyapounov's theorem and the Lynderberg-Levy-Feller's theorem. We give totally detailed proofs that are postponed to the (appendix) section \ref{sec_06}. In Section \ref{sec_03}, we deal with dependent random variables, here associated sequences. We make profit of the moving versions to significantly extend the central limit theorem of \cite{SangLo2018CLTA_01} which in turn is an extended version of \cite{GH-paulo}. Getting weak limits of \textit{MPS}'s requires the following more complicated condition

$$
\frac{1}{s_n^2} \mathbb{V}ar\left(\left\{\sum_{h=[ns_1]+1}^{[ns_2]} X_h\right\} + \left\{\sum_{h=[nt_1]+1}^{[nt_2]} X_h\right\}\right)\rightarrow \{a(s_2)-a(s_1)\}+\{a(t_2)-a(t_1)\},
$$ 

\Bin for $0=s_1<s_2<t_1<t_2\leq T \ (T=1)$, where $a(t)$ is a non-decreasing function in $t\in [0,1]$. Here again,  we apply the results in the weak limits in the finite-distribution invariance principle for associated sequences. In that section, we used regrouped-data method as it is usual in weak laws on associated sequences. 
Therein, we can use both conditions on regrouped data and no-regrouped data. We need a whole section, say Section \ref{sec_04}, to give the links between these two type of conditions. We close the paper with concluding remarks (in \ref{sec_05}).\\ 

\Ni The obtained results will help in successfully addressing the general setting of random invariance principles using random numbers $N(nt)$ of data, $t\in [0,1]$ in the innovative case of associated data, and beyond.\\

\Ni Let us proceed to the study for each type of dependence mentioned above.

\section{Central limit theorems for independent random variables} \label{sec_02}

\noindent In this section, we are going to check that Lyapounov's Theorem and Lynderberg-Levy-Feller's Theorem are unchanged in the moving frame. We exactly use the same proofs as in \cite{loeve} but we follow the detailed proofs in \cite{ips-mfpt-ang}. We will show how to use them in establishing  general invariance principles for independent data at least for finite-distributions. Here, we adopt the following notation
$$
s_{n}^{\prime2}=\sum_{k=p(n)+1}^{n+p(n)} \sigma_{k}^{2}, \ n\geq 1,
$$

\Bin and

$$
s_{n}^{2}=\sum_{k=1}^{n} \sigma_{k}^{2}, \ n\geq 1.
$$

\Bin Here are the moving versions of the two main central limit theorems for independent data.\\

\subsection{Moving versions}

\begin{theorem} \label{th_la} (Lyapounov). Suppose that the $X_k$'s are independent and have finite $(2+\delta)$-moments for every $\delta>0$ and 

\begin{equation} \label{cond_lya}
A_{n}^{\prime}(\delta) = \frac{1}{s_{n}^{\prime2+\delta}} \sum_{k=p(n)+1}^{n+p(n)}
\mathbb{E}\lvert X_{k}\rvert^{2+\delta}\longrightarrow 0, \ n\rightarrow +\infty.
\end{equation}

\bigskip \noindent Then, as $n\rightarrow +\infty$, 

$$
\frac{S_{n}^{\prime}}{s_{n}^{\prime}}\rightsquigarrow \mathcal{N}(0,1). \\
$$

\end{theorem}

\begin{remark}
\noindent As in the usual way, a version using arrays $\left\{X_{nk}, p(n)+1\leq k \leq n+p(n)\right\}$ is automatically written without any change in the proof.
\end{remark}

\begin{theorem} \label{th_LLF} (Feller-Levy-Lynderberg). Suppose that the $X_k$'s are independent and have only finite second order moments. We have the equivalence between the assertions below: \\ 

\begin{equation} \label{th_LLFa}
\underset{p(n)+1\leq k\leq n+p(n)}{\max} \left\{\frac{\sigma_{k}}{s_{n}^{\prime}}\right\}^{2}\longrightarrow0 \,\, and \,\, \frac{S_{n}^{\prime}}{s_{n}^{\prime}}\rightsquigarrow \mathcal{N}(0,1) , \ when \ n\rightarrow+\infty
\end{equation}

\bigskip \noindent and
  
\begin{equation} \label{th_LLFb}
\forall \epsilon>0, \  g_{n}(\epsilon) = \frac{1}{s_{n}^{\prime2}} \sum_{k=p(n)+1}^{n+p(n)} \int_{\left(|x|\geq \epsilon s_{n}^{\prime}\right)} x^{2} dF_{k}(x) \longrightarrow0, \ n\rightarrow +\infty. \\
\end{equation}
\end{theorem}

\Bin As promised, the proofs are direct and require checking all lines in the mentioned proofs. They are given in the appendix. Let us focus on the applications to invariance principles.\\

\Ni \textbf{Remark (R1)}. \label{remarkR1} A first remark is that if we have
$$
\limsup_{n\rightarrow +\infty} s_{n+p(n)}/s^{\prime}_{n}=:\nu \in ]0,+\infty[,
$$

\Bin then the Lynderberg condition for the whole sequence implies that it holds for the sequence $\{X_{p(n)+j}, j\geq 1\}$. Indeed, the condition above implies that for any $\eta>0$, there exists $n_0$ such that for any $n\geq n_0$, we have $s_{n+p(n)}<s^{\prime}_{n} (\nu+\eta)$. Hence, for $n\geq n_0$,

\begin{eqnarray*}
&&\frac{1}{s^{\prime 2}_{n}} \sum_{j=p(n)+1}^{n+p(n)} \int_{\left(|X_j|\geq \varepsilon s^{\prime}_{n}\right)} |X_j|^2 \ d\mathbb{P}\\
&&\leq \left(\frac{s^{2}_{n+p(n)}}{s^{\prime 2}_{n}}\right) \frac{1}{s^{2}_{n+p(n)}} \sum_{j=p(n)+1}^{n+p(n)} \int_{\left(|X_j|\geq \varepsilon^{\prime} s_{n+p(n)}\right)} |X_j|^2 \ d\mathbb{P}\\
&&\leq \left(\frac{s^{2}_{n+p(n)}}{s^{\prime 2}_{n}}\right) \frac{1}{s^{2}_{n+p(n)}} \sum_{j=1}^{n+p(n)} \int_{\left(|X_j|\geq \varepsilon^{\prime} s_{n+p(n)}\right)} |X_j|^2 \ d\mathbb{P},
\end{eqnarray*}

\Bin where $\varepsilon^{\prime}=\varepsilon (\nu + \eta)^{-1}$. By letting $n\rightarrow +\infty$, we get the moving Lynderberg condition. The same remark applies for the Lyapounov condition.\\

\subsection{Application to finite distributions limit invariance principles} $ $\\

\Ni Let us set \label{defSeqInvPrinc}

$$
\{Y_{n}(t), \ t \in [0,1]\}=:\left\{\frac{S_{[nt]}}{s_{n}} 1_{(nt\geq 1)}, \ t \in [0,1]\right\}.
$$

\Bin The invariance principle investigates whether such a sequence of stochastic processes converges to a tight stochastic process, mainly to a re-scaled Brownian motion. Here, a moving version of the central limit theorem is useful. It may be not very hard to proceed for independent data. But the way we use will serve as a basis for more complex dependent data as we will see in the second part of the paper.\\

\noindent The limit of $s_{[nt]}/s_n$ plays an important role here. In the \textit{iid} case with $\mathbb{E}X_1^2=\sigma^2$, we have $s_n^{2}=Var(S_n)=n \sigma^2$, which leads to the fact that

$$
(H0) \qquad \forall t \in \left[ 0, 1\right], \frac{s_{\left[nt\right]}^{2}}{s_n^{2}} \longrightarrow t, \qquad as  \qquad  n \longrightarrow +\infty. 
$$ 

\bigskip\noindent In the general case, we do not have such a simple relation. We have to set assumptions, for example
 
$$
(H1) \forall t \in \left[ 0, 1\right], \frac{s_{\left[nt\right]}^{2}}{s_n^{2}} \longrightarrow a(t), \qquad as \qquad n \longrightarrow +\infty,
$$ 

\Bin where $a(t)$ is a non-decreasing function of $t\in (0,1)$. From this, we may use the moving version of the Levy-Feller-Lynderberg's theorem to have the following result.

\begin{theorem}\label{invPrincIidFM} Suppose that the $X_j$'s are independent, centered and square integrable. Let us suppose that the Lynderberg condition holds. Then the sequence of stochastic processes $\{Y_n(t), 0\leq t \leq 1 \}$ weakly converges in finite distributions to 
$$
\left\{ W(a(t)),0\leq t\leq 1\right\},
$$

\Bin where $\{W(t), \ 0\leq t\leq 1\}$ is a Wiener stochastic process.
\end{theorem}

\Bin \textbf{Proof}. Let us set $0=t_0 <t_1 < t_2< ...< t_{k+1}=1$. We put : 

\begin{equation*}
\left\{ 
\begin{array}{c}
Z_{n}(t_{1})=Y_{n}(t_{1})-Y_{n}(t_{0})=\frac{1}{s_{n}}\sum_{[nt_{0}]<h\leq \lbrack nt_{1} \rbrack }X_{h} \\ 
\vdots \\ 
Z_{n}(t_{j})=Y_{n}(t_{j})-Y_{n}(t_{j-1})=\frac{1}{s_{n}}\sum_{[nt_{j-1}]<h\leq \lbrack nt_{j} \rbrack } X_{h} \\ 
\vdots\\ 
Z_{n}(t_{k})=Y_{n}(t_{k})-Y_{n}(t_{k-1})=\frac{1}{s_{n}}\sum_{[nt_{k-1}]<h\leq \lbrack nt_{k} \rbrack }X_{h}
\end{array}
\right..
\end{equation*}

\bigskip\noindent For each fixed $j \in \{1,...,k\}$, we have

$$
Z_{n}(t_{j})=\frac{1}{s_{n}}\sum_{h=[nt_{j-1}]+1}^{\lbrack nt_{j} \rbrack} X_{h}=\frac{\sqrt{s_{[nt_{j}]}^{2}-s_{[nt_{j-1}]}^{2}}}{\sqrt{s_{n}^{2}}}\times \frac{1}{\sqrt{s_{[nt_{j}]}^{2}-s_{[nt_{j-1}]}^{2}}} \sum_{h=[nt_{j-1}]+1}^{\lbrack nt_{j} \rbrack} X_{h}.
$$ 

\Bin Since the Lynderberg condition holds for the whole sequence, it holds for each sequence $\{X_{[nt_{j-1}]+h}, h\geq 1\}$, $j \in \{1,...,k\}$  [See Remark (R1), page 
\pageref{remarkR1}]. So, we may apply the moving Levy-Lynderberg-Feller's (CLT) for each $j \in \{1,\cdots,k\}$ with $p(n)=[nt_{j-1}]$ to get that

$$
\frac{1}{\sqrt{s_{[nt_{j}]}^{2}-s_{[nt_{j-1}]}^{2}}} \sum_{h=[nt_{j-1}]+1}^{\lbrack nt_{j} \rbrack} X_{h}
$$ 

\Bin weakly converges to the standard normal law $\mathcal{N}\left(0,1 \right)$ and by (H1) 

$$
\frac{s_{[nt_{j}]}^{2}}{s_{n}^{2}} \longrightarrow a(t_j) \ as \ n \longrightarrow +\infty.
$$  

\Bin So, for each $j \in \left\{1,\cdots, k \right\}$, by Slutsky's lemma, we have 

$$
Z_{n}(t_{j}) \rightsquigarrow \mathcal{N}\left(0, a(t_{j})-a(t_{j-1})\right), \ \ as \  n \longrightarrow +\infty,
$$ 

\bigskip\noindent since 

$$
\frac{s^{2}_{\left[nt_j\right]} -s^{2}_{\left[nt_{j-1}\right]}}{s^{2}_n}\rightarrow a(t_{j})-a(t_{j-1})
$$
	
\bigskip\noindent for each $j \in \{1,..., k\}$. Now, for any $u=(u_1, ..., u_k) \in \mathbb{R}^{k}$, 

\begin{equation*}
\mathbb{E}\left(\exp i\left(\sum_{1\leq j\leq k}Z_{n}(t_{j})u_{j}\right) \right)=\prod_{1\leq j\leq k}%
\mathbb{E}(\exp (iZ_{n}(t_{j})u_{j})\rightarrow \prod_{1\leq j\leq k}e^{-\frac{1}{2}u_{j}^{2}{\left(a(t_{j})-a(t_{j-1}) \right)}}.
\end{equation*}

\bigskip\noindent So, the vector $Z_n= \left( Z_n(t_j), 1\leq j \leq k  \right)^{t}$ weakly converges to a random Gaussian vector $Z$, with independent components. For each $1\leq j \leq k$ the $j^{th}$ component variance of $Z$ is: $a(t_j)-a(t_{j-1})$.\\

\noindent By the continuous theorem mapping theorem,  $Y_n=AZ_n$ with

\begin{equation*}
A=\left( 
\begin{array}{cccc}
1 & 0 & ... & 0 \\ 
1 & 1 & ... & 0 \\ 
1 & ... & 1 & 0 \\ 
1 & ... & ... & 1%
\end{array}
\right)
\end{equation*}

\bigskip\noindent weakly converges to $V=AZ$ and its components verify 

$$
V_j=Z_1+...+Z_j.
$$

\bigskip\noindent Next we see that 

$$
Z_j=V_j-V_{j-1}.
$$

\bigskip\noindent Hence, we conclude that the finite margins of  $Y_n(t)$ weakly converge to those of $W(a(t))$, where $W$ is a Wienner process, since

$$
\mathbb{E}(V_j^{2})=a(t_j) \qquad and  \qquad \mathbb{E}(V_j V_h)=a(t_j)\wedge a(t_h),
$$

\bigskip\noindent for all $j$ and $h$ both in $\{ 1,..,k \}$. \ $\blacksquare$\\


\section{central limit theorems for associated data} \label{sec_03}

\subsection{Easy introduction to Associated random variables} \label{sec_03_subsec_01}

\Ni In fear of rendering this paper heavier, we refer the reader to \cite{GH-gslo} for a quick introduction on associated sequence of random variables, and give the important Newman's inequality as follows. Readers are directed to \cite{rao} and \cite{bulinski} for a more detailed introduction to associations. The above mentioned inequality is used below several times.

\begin{lemma}[Newman and Wright (1981) Theorem, see \cite{GH-newmanwright}]
\label{GH-lemg3} Let $X_{1},X_{2},...,X_{n}$ be associated, then we have for
all $t=(t_{1},...,t_{n})\in \mathbb{R}^{n}$,

\begin{equation}
\left\vert \psi _{_{(X_{1},X_{2},...,X_{n})}}(t)-\prod\limits_{j=1}^{n}\psi
_{_{X_{j}}}(t_{j})\right\vert \leq \frac{1}{2}\sum_{1\leq j\neq h\leq
n}\left\vert t_{j}t_{h}\right\vert Cov(X_{j},X_{h}).  \label{GH-decomp}
\end{equation}
\end{lemma}

\Bin Let us begin with the approximation result.\\ 

\subsection{Moving version of the central limit theorem for associated data} \label{sec_03_subsec_02}$ $\\

\Ni As in \cite{SangLo2018CLTA_01}, we suppose that we have sequences of integer numbers $\ell (n)$,  $m(n)$, $r(n)$ such that
$n=m(n)\ell (n)+r(n),$ with $0\leq r(n)<\ell (n),$ $0\leq m(n)\rightarrow +\infty$ and

\begin{equation}
\ell (n)/n\rightarrow 0 \ as \ n\rightarrow +\infty .  \tag{L}
\end{equation}

\bigskip \noindent We want to stress that the integers $m=m(n),$ $\ell =\ell (n)$ and $r=r(n)$ are functions of $n\geq 1$ throughout the text even though we may and do drop the label $n$ in many situations for simplicity's sake. Let us suppose that $(p(n))_{n\geq 1}$ is an arbitrary sequence. The moving versions corresponding to the hypotheses in \cite{SangLo2018CLTA_01} and alike papers as \cite{GH-paulo} are as follows.\\

\Ni We denote

$$
S^{\prime}_0=0, \ \ S^{\prime}_n=\sum_{h=p(n)+1}^{n+p(n)} X_h \ \ and \ \ s^{\prime 2}_n=\mathbb{V}ar\left(S^{\prime}_n\right) \ for \ n\geq 1.
$$

\Bin The hypotheses we will be using are the following.\\

\begin{equation}
\frac{\ell (n)}{s^{\prime 2}_n}\rightarrow 0\text{ as }n\rightarrow +\infty; \tag{H0}
\end{equation}

\begin{equation}
\frac{\ell (n)}{s^{\prime 2}_{n}}\sum_{j=1}^{m(n)}\mathbb{V}ar\left( \frac{S^{\prime}_{j\ell(n)}-S^{\prime}_{(j-1)\ell (n)}}{\sqrt{\ell (n)}}\right) \rightarrow 1
\  as \ n\rightarrow +\infty;  \tag{Ha}
\end{equation}

\begin{equation}
\frac{1}{s^{\prime 2}_{n}}\mathbb{V}ar\left( \sum_{j=p(n)+m(n)\ell(n)+1}^{p(n)+n} X_j\right) \rightarrow 0\text{ as }n\rightarrow +\infty;
\tag{Hab}
\end{equation}

\noindent 

\begin{equation}
\sup_{1\leq j\leq m(n)}\frac{\ell (n)}{s^{\prime 2}_{n}}\mathbb{V}ar\left(\frac{S^{\prime}_{j\ell (n)}-S^{\prime}_{(j-1)\ell (n)}}{\sqrt{\ell (n)}}\right)
=C_{1}(n)\rightarrow 0\text{ as \ }n\rightarrow +\infty .  \tag{Hb}
\end{equation}

\Bin
\noindent In the sequel, it may be handy to use the notation 
\begin{equation}
Y_{j,\ell }=\frac{S^{\prime}_{j\ell (n)}-S^{\prime}_{(j-1)\ell (n)}}{\sqrt{\ell (n)}},1\leq
j\leq m=m(n).  \label{GH-THEY}
\end{equation}

\Bin Let us prove the moving version of Formula (4.12) in \cite{SangLo2018CLTA_01}.

\begin{proposition} \label{movingApproxAsso_independt_01}
Under Hypotheses (L), (H0), (Ha), (Hab) and (Hb), we have for any $t\in \mathbb{R}$, as $n\rightarrow +\infty$,
\begin{equation}
\left\vert \Psi _{\frac{S^{\prime}_{n}}{s^{\prime}_{n}}}(t)-\prod\limits_{j=1}^{m}\Psi_{Y_{j,\ell}}\left( \frac{\sqrt{\ell }}{s^{\prime}_{n}}t\right) \right\vert \rightarrow 0.  \label{GH-conclusionFin}
\end{equation}
\end{proposition}

\Bin \textbf{Proof}. We follow the one in \cite{SangLo2018CLTA_01} which we appropriately adapt. Let us denote

$$
\Psi _{\frac{S^{\prime}_{n}}{s^{\prime}_{n}}}(t)=\mathbb{E}\left( e^{itS^{\prime}_{n}/s^{\prime}_{n}}\right), t\in \mathbb{R}.
$$

\Bin We have for $t\in \mathbb{R}$, 
\begin{eqnarray}
\left\vert \Psi _{_{\frac{S^{\prime}_{n}}{s^{\prime}_{n}}}}(t)-\Psi _{_{\frac{S^{\prime}_{m\ell }}{s^{\prime}_{n}}}}(t)\right\vert &=&\left\vert \mathbb{E}(e^{itS^{\prime}_{n}/s^{\prime}_{n}})-\mathbb{E}(e^{itS^{\prime}_{m\ell }/s^{\prime}_{n}})\right\vert\\
&=&\left\vert \mathbb{E}\left[ e^{itS^{\prime}_{m\ell }/s_n^{\prime}}\left( e^{it\left[
(S^{\prime}_{n}/s^{\prime}_{n})-(S^{\prime}_{m\ell }/s^{\prime}_{n})\right] }-1\right) \right] \right\vert\\
&\leq& \mathbb{E}\left\vert e^{it\left( \frac{S^{\prime}_{n}}{s^{\prime}_{n}}-\frac{S^{\prime}_{m\ell }}{
s^{\prime}_{n}}\right) }-1\right\vert .  \label{GH-b}
\end{eqnarray}

\noindent But for any $x\in \mathbb{R}$, 
\begin{equation*}
\left\vert e^{ix}-1\right\vert =|(\cos x-1)+i\sin x|=\left\vert 2\sin \frac{x%
}{2}\right\vert \leq |x|.
\end{equation*}

\bigskip \noindent Thus the second member of $(\ref{GH-b})$ is, by the
Cauchy-Schwarz's inequality, bounded by 
\begin{equation*}
|t|\mathbb{E}\left\vert \frac{S^{\prime}_{n}}{s^{\prime}_{n}}-\frac{S^{\prime}_{m\ell }}{s^{\prime}_{n}}%
\right\vert \leq |t|\mathbb{V}ar\left( \frac{S^{\prime}_{n}}{s^{\prime}_{n}}-\frac{S^{\prime}_{m\ell }}{s^{\prime}_{n}}\right)^{\frac{1}{2}}
\end{equation*}

\noindent and 
\begin{equation*}
\delta _{m,\ell }=\mathbb{V}ar\left( \frac{S^{\prime}_{n}}{s^{\prime}_{n}}-\frac{S^{\prime}_{m\ell }}{s^{\prime}_{n}}\right) =\frac{1}{s^{\prime 2}_{n}}\mathbb{V}ar\left(S^{\prime}_{n}-S^{\prime}_{m\ell}\right),
\end{equation*}

\Bin which tends to zero as $n\rightarrow +\infty $  by \textit{(Hab)} since%

\begin{eqnarray*}
\delta _{m,\ell } &=&\frac{1}{s^{\prime 2}_{n}}\mathbb{V}ar\left(S^{\prime}_{n}-S^{\prime}_{m\ell}\right)   \label{GH-C1A} \\
&= &\frac{1}{s^{\prime 2}_{n}}\mathbb{V}ar\left( \sum_{j=p(n)+m(n)\ell(n)+1}^{p(n)+n} X_j\right)  \\
& \rightarrow & 0.
\end{eqnarray*}

\Bin  This proves that
 
\begin{equation}
\left|\Psi _{\frac{S^{\prime}_{n}}{s^{\prime}_{n}}}(t)-\Psi _{\frac{S^{\prime}_{m\ell }}{s^{\prime}_{n}}%
}(t)\right|\rightarrow 0\text{ as }n\rightarrow +\infty .  \label{GH-commonStep01}
\end{equation}

\Bin Next, recall that  $Y_{j,\ell }=(S^{\prime}_{j\ell }-S^{\prime}_{(j-1)\ell})/\sqrt{\ell }$%
, for $1\leq j\leq m$. \ Observe that

\begin{equation*}
\frac{S^{\prime}_{m\ell }}{s^{\prime}_{n}}=\frac{\sqrt{\ell }}{s^{\prime}_{n}}\sum_{j=1}^{m}Y_{j,\ell }.
\end{equation*}

\Bin According to the Newman's inequality (see Lemma \ref{GH-lemg3}), we have
\begin{equation*}
\left\vert \Psi _{\frac{S^{\prime}_{m\ell }}{s^{\prime}_{n}}}(t)-\prod\limits_{j=1}^{m}\Psi
_{Y_{j,\ell}}\left( \frac{\sqrt{\ell }}{s^{\prime}_{n}}t\right) \right\vert \leq \frac{%
\ell t^{2}}{2s^{\prime 2}_{n}}\sum_{1\leq j\neq k\leq m}Cov(Y_{j,\ell },Y_{k,\ell
}).
\end{equation*}

\noindent But, 
\begin{eqnarray*}
\frac{\ell t^{2}}{2s^{\prime 2}_{n}}\sum_{1\leq j\neq k\leq m}Cov(Y_{j,\ell},Y_{k,\ell }) &=&\frac{\ell t^{2}}{2s^{\prime 2}_{n}}\mathbb{V}ar\left(\sum_{j=1}^{m}Y_{j,\ell }\right) -\frac{\ell t^{2}}{2s^{\prime 2}_{n}}\sum_{j=1}^{m}\mathbb{V}ar(Y_{j,\ell }) \\
&=&\frac{t^{2}}{2}\left[ \mathbb{V}ar\left( \frac{\sqrt{\ell }}{s^{\prime}_{n}}\sum_{j=1}^{m}Y_{j,\ell }\right) -\frac{\ell }{s^{\prime 2}_{n}}\sum_{j=1}^{m}\mathbb{V}ar\left( Y_{j,\ell }\right) \right]  \\
&=&\frac{t^{2}}{2}\left[ \mathbb{V}ar\left( \frac{1}{s^{\prime}_{n}}S^{\prime}_{m\ell }\right)-\frac{\ell }{s^{\prime 2}_{n}}\sum_{j=1}^{m}\mathbb{V}ar\left( \frac{S_{j\ell}-S_{(j-1)\ell}}{\sqrt{\ell }}\right) \right].  (L3)\\
\end{eqnarray*}

\Bin From (L3), we use

$$
\frac{1}{s^{\prime}_{n}}S^{\prime}_{n}= \frac{1}{s^{\prime}_{n}}S^{\prime}_{m\ell} +\frac{1}{s^{\prime}_{n}} \sum_{j=p(n)+m\ell+1}^{n+p(n)}X_{j},
$$

\Bin and take the associativity into account to get

$$
1=\mathbb{V}ar\left(\frac{1}{s^{\prime}_{n}}S^{\prime}_{n}\right) \geq \mathbb{V}ar\left(\frac{1}{s^{\prime}_{n}}S^{\prime}_{m\ell}\right) + \mathbb{V}ar\left(\frac{1}{s^{\prime}_{n}} \sum_{j=p(n)+m\ell+1}^{n+p(n)}X_{j}\right).
$$

\Bin This leads to

\begin{eqnarray*}
&&\frac{\ell t^{2}}{2s^{\prime 2}_{n}}\sum_{1\leq j\neq k\leq m}Cov(Y_{j,\ell},Y_{k,\ell })\\
&&\leq\frac{t^{2}}{2}\left[ 1-\frac{\ell }{s^{\prime 2}_{n}}\sum_{j=1}^{m}\mathbb{V}ar\left( \frac{S^{\prime}_{j\ell }-S^{\prime}_{\ell (j-1)}}{\sqrt{\ell }}\right) \right]\\
&&-\frac{t^{2}}{2s^{\prime 2}_{n}}\mathbb{V}ar\left(\sum_{j=1}^{r}X_{p(n)+m\ell+j}\right).
\end{eqnarray*}

\Bin Thus, by \textit{(Ha)} and \textit{(Hab)}, we get, as  $n\rightarrow +\infty$ our conclusion, i.e., 

\begin{equation}
\left\vert \Psi _{\frac{S^{\prime}_{m\ell }}{s^{\prime}_{n}}}(t)-\prod\limits_{j=1}^{m}\Psi
_{Y_{j,\ell}}\left( \frac{\sqrt{\ell }}{s^{\prime}_{n}}t\right) \right\vert \rightarrow 0 \ as \ n\rightarrow +\infty .  \label{GH-conclusion01}
\end{equation}

\Bin Let us draw a first important application of our results regarding the invariance principle.\\

\subsection{Gaussian Central limit theorem} \label{sec_03_subsec_03}

\Ni Proposition \ref{movingApproxAsso_independt_01} says that under the hypotheses (L), (H0), (Ha), (Hab) and (Hb), $S^{\prime}_n/s^{\prime}_n$ has the same weak limit law than

$$
\frac{1}{s^{\prime}_n}\sum_{j=1}^{m(n)} T_{j,n},
$$

\Bin where, for each $n\geq 1$, $(T_{1,n}, \cdots, T_{m(n),n})$ has independent components and for each $1\leq j \leq m(n)$

$$
T_{j,n}=_{d} S^{\prime}_{j\ell (n)}-S^{\prime}_{(j-1)\ell (n)},
$$

\Bin where the random variables in the left-hand of the equality are already defined in Formula \eqref{GH-THEY} and $=_{d}$ stands for the equality in distribution.\\

\Ni \textbf{A - Lyapounov central limit theorem}. As a consequence a Lyapounov condition for independent data is enough to ensure the Gaussian central limit theorem. Let us set for $\delta>0$,

\begin{equation}
\frac{1}{s^{\prime (2+\delta)}_{n}}\sum_{j=1}^{m}\mathbb{E}\left\vert S^{\prime}_{j\ell (n)}-S^{\prime}_{(j-1)\ell (n)}\right\vert^{2+\delta }=C_{2}(n)\rightarrow 0\text{ as n}\rightarrow +\infty .  \tag{Hc}
\end{equation}

\Bin By the $C_{r}$-inequality ($|a+b|^r\leq 2^{r-1}(|a|^{r}+|b|^{r})$ for real-valued numbers $a$, $b$ and $r\geq 1$), the $S^{\prime}_{j\ell (n)}-S^{\prime}_{(j-1)\ell (n)}$'s have finite $(2+\delta)$-moments, $\delta>0$ if the $X_j$'s do.\\

\Bin So, we have

\begin{theorem} \label{assoc_lyapounov} Let $\delta>0$. If the  $X_j$'s have finite $(2+\delta)$-moments and (Hc) holds on top of the hypotheses (L), (H0), (Ha), (Hab) and (Hb), then

$$
S^{\prime}_n/s^{\prime}_n \rightsquigarrow \mathcal{N}(0,1).
$$
\end{theorem}

\Bin \textbf{Proof} Since the study transformed into that of some in independent and centered data, the Lyapounov condition is enough to get the conclusion, but the involved variance is $s^{\prime 2}_{m\ell}$ which is equivalent to $s^{\prime 2}_n$ by Hypotheses (Ha).\\


\Ni \textbf{B - A Feller-Levy-Lynderberg central limit theorem}.\\

\Ni With the equivalence $s^{\prime 2}_{m\ell} \sim s^{\prime 2}_n$ by Hypothesis (Ha), the Lynderberg condition on the $T_{j,n}$'s is equivalent to: $\forall$ $\epsilon>0$, 

\begin{equation}
\frac{1}{s^{\prime 2}_n} \sum_{j=1}^{m(n)} \int_{(|S^{\prime}_{j\ell (n)}-S^{\prime}_{(j-1)\ell (n)}|\geq \varepsilon \tau^{\prime}_n)} |S^{\prime}_{j\ell (n)}-S^{\prime}_{(j-1)\ell (n)}|^2 \ \ d\mathbb{P}\longrightarrow 0, \label{ynderAssocA}
\end{equation}

\Bin where

$$
\tau^{\prime 2}_{j,n}=\mathbb{V}ar(S^{\prime}_{j\ell (n)}-S^{\prime}_{(j-1)\ell (n)}), \ 1\leq j \leq m(n)
$$

\Bin and

$$
\tau^{\prime 2}_{n}=\sum_{j=1}^{m(n)} \tau^{\prime 2}_{j,n}.
$$

\Bin We denote

$$
B_n= \max\{\tau^{\prime}_{j,n}/\tau^{\prime}_{n}, \ 1 \leq j \leq m(n)\}.
$$

\Bin The Feller-Levy-Lynderberg (FLL) theorem is stated as fallows.\\

\begin{theorem} \label{assoc_FLL} We suppose that the $X_j$'s have finite second order moments and that the hypotheses (L), (H0), (Ha), (Hab) and (Hb) hold. Then

$$
S^{\prime}_n/s^{\prime}_n \rightsquigarrow \mathcal{N}(0,1) \ \ and \ \ B_n\rightsquigarrow 0
$$

\Bin if and only if the Lynderberg condition \eqref{ynderAssocA} holds.

\end{theorem}

\Bin Here, we directly used the Lyapounov and the FLL conditions ((Hc) and  \ref{ynderAssocA}) for on the regrouped data (the $T_{j,n}$'s). It may be more convenient to give sufficient conditions of the non-regrouped data for they hols. When we proceed to complete comparison between the situation for regrouped data and non-regrouped data as we will do it in Section \label{sec_04}, we arrive at this final version.\\

\begin{theorem} \label{general_WL_MPS_ASSOC} Suppose that the random variables $X_j$, $j\geq 1$, are centered and square integrable and associated and that Hypotheses (L), (H0), (Ha), (Hab) and (Hb) are satisfied. Then we have following results:\\

\Ni (1) If the $X_j$'s have $(2+\delta)$-moments for some $\delta>0$ and the Lyapounov-type condition 

\begin{equation}
\frac{\ell(n)^{1+\delta}}{s^{\prime (2+\delta)}_n} \sum_{p(n)+1 \leq j \leq p(n)+n} \mathbb{E} \left|X_{j}\right|^{2+\delta} \rightarrow 0 \ as \ n\rightarrow +\infty, \label{L-NRD}
\end{equation}

\Bin then 

$$
S^{\prime}_n/s^{\prime}_n \rightsquigarrow \mathcal{N}(0,1) \ as \ n\rightarrow +\infty.
$$

\Bin (2) Suppose that the following uniform negligibility of the variances

$$
\frac{\ell(n)^2}{\epsilon^2 s^{\prime 2}_n} \max_{p(n)+1 \leq j \leq p(n)+n}  \mathbb{E}X_{j}^2 \rightarrow 0 \ as \ n\rightarrow +\infty.
$$

\Bin hold. Then we have

$$
S^{\prime}_n/s^{\prime}_n \rightsquigarrow \mathcal{N}(0,1) \ as \ n\rightarrow +\infty.
$$

\Bin if and only if the following Lynderberg-type condition holds:

\begin{equation}
L_n^{\prime}\left(\frac{\varepsilon}{2\ell(n)}\right) \rightarrow 0 \ as \ n\rightarrow +\infty. \label{LYND-NRD}
\end{equation}

\Bin where, for $\varepsilon>0$ and $n \geq 1$, 

\begin{equation*}
L^{\prime}_n(\varepsilon)=\frac{\ell(n)^2}{s_n^{\prime 2}} \sum_{p(n)+1 \leq j \leq p(n)+n} \int_{(|X_j| >\varepsilon s^{\prime}_n)} X_j^2 \ d\mathbb{P}.
\end{equation*}
\end{theorem}

\Bin As said earlier, we will give a complete justification of that theorem in Section \ref{sec_04}.\\

\Ni Now, let us see the applications of the results to invariance principles.\\

\subsection{Invariance principles of associated data} \label{sec_03_subsec_03}

\Ni We already defined the sequence of stochastic processes 

$$
\{Y_{n}(t), \ t \in [0,1]\}=:\left\{\frac{S_{[nt]}}{s_{n}} 1_{(nt\geq 1)}, \ t \in [0,1]\right\}
$$

\Bin in page \pageref{defSeqInvPrinc} and we wish to find the weak law. Above, in the independent case, we used the hypothesis (H1) which controls how the variance $s_{[nt]}$ grows with respect to $s_{n}$. Such a condition still works if the data are stationary. Otherwise, we may need the following more elaborated assumption.\\

\Ni \textbf{(H1-NSA)}. There exists a measurable function $a(t)$ of $t \in [0,1]$ such that for $0<s_1<s_2\leq t_1<t_2\leq 1$, as $n\rightarrow +\infty$,

$$
\frac{1}{s_n^2} \mathbb{V}ar\left(\left\{\sum_{h=[ns_1]+1}^{[ns_2]} X_h\right\} + \left\{\sum_{h=[nt_1]+1}^{[nt_2]} X_h\right\}\right)\rightarrow \{a(s_2)-a(s_1)\}+\{a(t_2)-a(t_1)\}. 
$$ 

\Bin From this, we may use the moving version of the Levy-Feller-Lynderberg's theorem or Lyapounov theorem, we have the following result.

\begin{theorem}\label{invPrincINSA} Let us suppose that the $X_j$'s are associated and (H1-NSA) holds. Let us suppose that the hypotheses (L), (H0), (Ha), (Hab) and (Hb) hold on top of the moving Lynderberg condition. Then the sequence of stochastic processes $\{Y_n(t), 0\leq t \leq 1 \}$ weakly converges in finite distributions to 

$$
\left\{ W(a(t)),0\leq t\leq 1\right\},
$$

\Bin where $\{W(t), \ 0\leq t\leq 1\}$ is a Wienner stochastic process.
\end{theorem}

\noindent \textbf{Proof}. Let us set $0=t_0 <t_1 < t_2< ...< t_{k+1}=1$. We put : 

\begin{equation*}
\left\{ 
\begin{array}{c}
Z_{1,n}=Y_{n}(t_{1})-Y_{n}(t_{0})=\frac{1}{s_{n}}\sum_{[nt_{0}]<h\leq \lbrack nt_{1} \rbrack }X_{h} \\ 
\vdots \\ 
Z_{j,n}=Y_{n}(t_{j})-Y_{n}(t_{j-1})=\frac{1}{s_{n}}\sum_{[nt_{j-1}]<h\leq \lbrack nt_{j} \rbrack } X_{h} \\ 
\vdots\\ 
Z_{k,n}=Y_{n}(t_{k})-Y_{n}(t_{k-1})=\frac{1}{s_{n}}\sum_{[nt_{k-1}]<h\leq \lbrack nt_{k} \rbrack }X_{h}%
\end{array}
\right. 
\end{equation*}

\bigskip\noindent For each fixed $j \in \{1,...,k\}$, $s_n Z_{j,n}$ is a moving partial sum with $p(j,n)=[nt_{j-1}]$, $\Delta(j,n)=[nt_{j}]-[nt_{j-1}]$ and

$$
s_n Z_{j,n}=\sum_{h=p(j,n)+1}^{p(j,n)+\Delta(j,n)} X_h.
$$

\Bin Put, for each $j \in \{1,...,k\}$,  

$$
s^{\prime 2}_{j,n}=\mathbb{V}ar\left(\sum_{h=p(j,n)+1}^{p(j,n)+\Delta(j,n)} X_h\right).
$$

\Bin We have the central limit theorem for $s_n Z_{j,n}$, i.e.,

$$
\frac{s_n Z_{j,n}}{s^{\prime}_{j,n}} \rightsquigarrow \mathcal{N}(0,1),
$$

\Bin and by Slutsky's lemma, we have 

$$
Z_{j,n} \rightsquigarrow \mathcal{N}(0,\Delta a(t_j)),
$$

\Bin where we denote $\Delta a(t_j)=a(t_j)-a(t_{j-1})$, $j \in \{1,...,k\}$. Next, let us apply Wold Criterion (see \cite{ips-mfpt-ang}, Chapter 1) and consider

$$
Z_n=u_1 Z_{1,n} + \cdots + u_k Z_{k,n}.
$$

\Bin For  $t\in \mathbb{R}$, we have

\begin{eqnarray*}
&&\left|\psi_{Z_n}(t) -\prod_{j=1}^{k} \exp\left(-\frac{1}{2}\Delta a(t_j)t^2 u_j^2\right) \right|\\
&&\leq \left|\psi_{Z_n}(t) -\prod_{j=1}^{k} \psi_{Z_{j,n}}(tu_j) \right|\\
&&+\left|\prod_{j=1}^{k} \psi_{Z_{j,n}}(tu_j) -\prod_{j=1}^{k} \exp\left(-\frac{1}{2}\Delta a(t_j)t^2 u_j^2\right) \right|\\
&& =: R_n(1) + R_n(2).
\end{eqnarray*}

\Bin We already proved that $R_n(2) \rightarrow 0$. It remains to prove that $R_n(1)$ does. But, by Newman's inequality in Lemma \ref{GH-lemg3}, we have

$$
0\leq R_n(1) \leq \frac{t^2}{4} \sum_{1\leq j_1<j_2\leq k} |u_{j_1}u_{j_2}| \mathbb{C}ov\left(Z_{j_1,n}, Z_{j_2,n}\right).
$$

\Bin Now for any $1\leq j_1<j_2\leq k$, 

$$
2\mathbb{C}ov\left(Z_{j_1,n}, Z_{j_2,n}\right)=R_n(1,1)-R_n(1,2),
$$

\Bin with

$$
R_n(1,1)=\frac{1}{s^2_n} \left(\mathbb{V}ar\left(\left\{\sum_{h=[nt_{j_1-1}]+1}^{[nt_{j_1}]} X_h\right\} + \left\{\sum_{h=[nt_{j_2-1}]+1}^{[nt_{j_2}]} X_h\right\}\right)\right)
$$

\Bin and 

$$
R_n(1,2)=\frac{1}{s^2_n} \mathbb{V}ar\left(\sum_{h=[nt_{j_1-1}]+1}^{[nt_{j_1}]} X_h\right) + \frac{1}{s^2_n} \mathbb{V}ar\left(\sum_{h=[nt_{j_2-1}]+1}^{[nt_{j_2}]} X_h\right),
$$

\Bin which both converge to $\Delta a(t_{j_1}) + \Delta a(t_{j_2})$ and hence $2\mathbb{C}ov\left(Z_{j_1,n}, Z_{j_2,n}\right)\rightarrow 0$. From there, the conclusion is the same as in the independence case since we move from the limit to the increments to the sequence $(Y_n(t_1), \cdots, Y_n(t_k))$ in the same way. $\square$\\

\newpage

\section{Weak limits of moving partial sums using condition on full data} \label{sec_04}

\Ni Theorems \ref{assoc_lyapounov} and \ref{assoc_FLL} used Lyapounov and \textit{FLL}-type conditions on the regrouped date. Using Proposition \ref{movingApproxAsso_independt_01}, we pointed out that  under Hypotheses (H0), (Ha), (Hab) and (Hb), the weak limit law of $Y_n=S^{\prime}_n/s^{\prime}_n$ behaves exactly as that of

$$
T^{\prime\prime}_n=\sum_{j=1}^{m(n)} T_{j,n}/\tau^{\prime}_n,
$$ 

\Bin where $T_{j,n}=_d S^{\prime}_{j\ell(n)}-S^{\prime}_{(j-1)\ell(n)}$, $1\leq j \leq m(n)$, and the $T_{j,n}$ are independent, with 

$$
\tau^{\prime 2}_n=\sum_{j=1}^{m(n)} \tau^{\prime 2}_{j,n}, \ \ \tau^{\prime 2}_{j,n}=\mathbb{V}ar(T_{j,n}), \ 1\leq j\leq m(n). 
$$

\Ni From there, the theory on weak limits of independent case applies. In studying $T^{\prime\prime}_n$ we have the Bounded Variance Hypothesis (BVH) is satisfied with

$$
\sup_{n\geq 1} \sum_{1\leq j\leq m(n)} \mathbb{V}ar(T_{j,n}/\tau^{\prime}_n)\leq c,
$$

\Bin for $c=1$. But, because of Hypothesis (Ha), we may change the sequence we study to

$$
T^{\prime\prime}_n=\sum_{j=1}^{m(n)} T_{j,n}/s^{\prime}_n
$$

\Bin and the \textit{(\textit{BVH})} still holds for some $c>0$. Moreover, the Uniformly Asymptotically Negligibility (\textit{UAN}) condition in sums of independent random variables theory :

\begin{equation}
U^{\prime\prime}(n)= \max_{1\leq j \leq m(n)} \mathbb{P}(|T_{j,n}|\geq \varepsilon \tau^{\prime}_n) \rightarrow 0,
\end{equation}

\Bin is controlled as follows:

\begin{equation}
\forall \varepsilon>0, \ \ U^{\prime\prime}(n)\leq B^{\prime\prime}_n =\frac{\varepsilon^{-2}}{\tau^{\prime 2}_{n}}\max_{\ 1\leq j \leq m(n)} \tau^{\prime 2}_{j,n} \rightarrow 0. \label{BN-R}
\end{equation}

\Bin Now, in the first place, let us see the links of usual conditions related to theory of sums of independent random variables when they are expressed on the $T_{j,n}$'s and on the $X_j$'s.\\

\Ni \textbf{A - Convergence conditions on the regrouped and simple data}.\\

\Ni \textbf{(A1) The \textbf{UAN} condition}. That condition is expressed on the regrouped data $T_{j,n}$'s. It is controlled with $B^{\prime\prime}_n$. This itself is bounded through the simple data $\{X_{p(n)+j}, \ 1\leq j \leq n\}$ as follows. By using the convexity of $\mathbb{R}_+ \ni x \rightarrow x^{2}$, we have

\begin{eqnarray*}
B^{\prime\prime}_n&=&\frac{1}{\epsilon^2\tau^{\prime 2}_n}\max_{1\leq j \leq m(n)} \mathbb{E} \left(\sum_{h=1}^{\ell(n)} X_{p(n)+(j-1)\ell(n)+h}\right)^2\\
&\leq&\frac{\ell(n)}{\epsilon^2\tau^{\prime 2}_n}\max_{1\leq j \leq m(n)}  \sum_{h=1}^{\ell(n)} \mathbb{E}X_{p(n)+(j-1)\ell(n)+h}^2\\
&\leq&\frac{\ell(n)^2}{\epsilon^2\tau^{\prime 2}_n} \max_{p(n)+1 \leq j \leq p(n)+n}  \mathbb{E}X_{j}^2.\\
\end{eqnarray*}

\Bin By denoting 

\begin{equation}
B_n^{\prime}=\frac{\ell(n)^2}{\epsilon^2 s^{\prime 2}_n} \max_{p(n)+1 \leq j \leq p(n)+n}  \mathbb{E}X_{j}^2, 
\label{BN-S}
\end{equation}

\Bin we have

\begin{equation}
B_n^{\prime\prime} \leq\left(\frac{s^{\prime}_n}{\tau^{\prime}_n}\right)^2 B_n^{\prime}. \label{BN-RS}
\end{equation}

\Bin \textbf{(A2) - Lyapounov condition}. For $\delta>0$, the Lyapounov condition is

\begin{equation}
A^{\prime\prime}_n(\delta)=\frac{1}{\tau^{\prime (2+\delta)}_n} \sum_{1\leq j \leq m(n)} \mathbb{E} \left|\sum_{1\leq h\leq \ell(n)} X_{p(n)+(j-1)\ell(n)+h}\right|^{2+\delta}. \label{Lyap-R}
\end{equation}

\Bin To shorten the notation, we use the following notation below

$$
X^{\prime}_{n,j,h}=X_{p(n)+(j-1)\ell(n)+h}, \ 1\leq j\leq m(n), \ \ 1\leq h \leq \ell(n).
$$

\Bin Now, by using the convexity of $\mathbb{R}_+ \ni x \rightarrow x^{2+\delta}$, we have

\begin{eqnarray*}
A^{\prime\prime}_n(\delta)&=&\frac{1}{\tau^{\prime (2+\delta)}_n} \sum_{1\leq j \leq m(n)} \mathbb{E} 
\left| \ell(n) \left( \sum_{1\leq h\leq \ell(n)} \{ X^{\prime}_{n,j,h}/\ell(n)\}\right)\right|^{2+\delta}\\
&\leq& \frac{\ell(n)^{1+\delta}}{\tau^{\prime (2+\delta)}_n} \sum_{1\leq j \leq m(n)} \sum_{1\leq h\leq \ell(n)} \mathbb{E} \left|X^{\prime}_{n,j,h}\right|^{2+\delta}\\
&\leq& \frac{\ell(n)^{1+\delta}}{\tau^{\prime (2+\delta)}_n} \sum_{p(n)+1 \leq j \leq p(n)+n} \mathbb{E} \left|X_{j}\right|^{2+\delta},\\
\end{eqnarray*}

\Bin and we get

\begin{equation}
A^{\prime\prime}_n(\delta) \leq  \left(\frac{s^{\prime (2+\delta)}_n}{\tau^{\prime (2+\delta)}_n}\right) A^{\prime}_n(\delta), \label{Lyap-RS}
\end{equation}

\Bin where

\begin{equation}
A^{\prime}_n(\delta)=
\frac{\ell(n)^{1+\delta}}{s^{\prime (2+\delta)}_n} \sum_{p(n)+1 \leq j \leq p(n)+n} \mathbb{E} \left|X_{j}\right|^{2+\delta}. \label{Lyap-S}
\end{equation}

\Bin \textbf{(A3) - Lynderberg Condition}. For $\varepsilon>0$, the Lynderberg Condition for regrouped data is

\begin{equation}
L^{\prime\prime}_n(\varepsilon)=\frac{1}{\tau^{\prime 2}_n} \sum_{j=1}^{m(n)} \int_{\left(\left|T_{j,n}\right|\geq \varepsilon \tau^{\prime}_n\right)} \left|T_{j,n}\right|^2 \ d\mathbb{P}. \label{Lynder-R}
\end{equation}

\Bin Let us set, $n\geq 1$, $1\leq j \leq m(n)$, $1\leq h \leq \ell(n)$,

$$
A^{\prime}_{n,j,h}=\left(\left|X^{\prime}_{n,j,h}\right|\geq \varepsilon \tau^{\prime}_n/\ell(n)\right)
$$

\Bin and

$$
A^{\prime}_{n,j}= \bigcup_{j=1}^{\ell(n)} A^{\prime}_{n,j,h}.
$$

\Bin We have

\begin{equation*}
L^{\prime\prime}_n(\varepsilon)\leq \frac{1}{\tau^{\prime 2}_n} \sum_{j=1}^{m(n)} \int_{A^{\prime}_{n,j}} \left|T_{j,n}\right|^2 \ d\mathbb{P}. 
\end{equation*}

\Bin Let $M_n=\max_{1\leq s\leq \ell(n)} |X^{\prime}_{n,j,s}|$. Since we have

$$
\bigcup_{1\leq r \leq \ell(n)} \biggr(M_n=|X^{\prime}_{n,j,r}|\biggr)=\Omega,
$$

\Bin we get

\begin{eqnarray*}
L^{\prime\prime}_n(\varepsilon)&\leq& \frac{1}{\tau^{\prime 2}_n} \sum_{j=1}^{m(n)} \sum_{r=1}^{\ell(n)} \int_{A^{\prime}_{n,j} \cap (M_n=|X^{\prime}_{n,j,r}|)} \left|T_{j,n}\right|^2 \ d\mathbb{P}\\
&\leq& \frac{1}{\tau^{\prime 2}_n} \sum_{j=1}^{m(n)} \sum_{r=1}^{\ell(n)} \int_{A^{\prime}_{n,j} \cap (M_n=|X^{\prime}_{n,j,r}|)} (\ell(n) |X^{\prime 2}_{n,j,r}|)^2 \ d\mathbb{P}\\
&\leq& \frac{\ell(n)^2}{\tau^{\prime 2}_n} \sum_{j=1}^{m(n)} \sum_{r=1}^{\ell(n)} \int_{(|X^{\prime}_{n,j,r}|> \varepsilon \tau^{\prime}_n/\ell(n))} X^{\prime 2}_{n,j,r} \ d\mathbb{P}\\
&\leq& \frac{\ell(n)^2}{\tau^{\prime 2}_n} \sum_{j=p(n)+1}^{p(n)+n}  \int_{(|X_j|> \varepsilon \tau^{\prime}_n/\ell(n))} X_j^2 \ d\mathbb{P}\\
\end{eqnarray*}

\Bin We are going to complete our comparison by using a moving Lynderberg condition on the data $X_j$'s:

\begin{equation}
L^{\prime}_n(\varepsilon)=\frac{\ell(n)^2}{s_n^{\prime 2}} \sum_{p(n)+1 \leq j \leq p(n)+n} \int_{(|X_j| >\varepsilon s^{\prime}_n)} X_j^2 \ d\mathbb{P}. \label{Lynder-O-1}
\end{equation}

\Bin But, by Hypothesis (\textit{Ha}) and (\textit{Hab}), we have  $(\tau_n^{\prime})^2/s_n^{\prime 2}\rightarrow 1$, as $n\rightarrow +\infty$. So, for $n$ large enough, $(|X_j| >\varepsilon \tau^{\prime}_n)$ implies $(|X_j| >(\varepsilon/2) s^{\prime}_n)$. So, we have for $n$ large enough, for any $\varepsilon>0$,

\begin{equation}
L^{\prime\prime}_n(\varepsilon) \leq \left(\frac{s^{\prime}_n}{\tau^{\prime}_n}\right)^2 L^{\prime}_n(\varepsilon/(2\ell(n))). \label{Lynder-O-2}
\end{equation}

\Bin Now we may justify Theorem \ref{general_WL_MPS_ASSOC} using conditions of the real data $X_j$, $j\geq 1$, as announced, by summarizing the discussion above.\\

\Bin \textbf{(B) Weak limits of Moving Partial sums using no-regrouped data}.\\

\Ni Let us suppose that hypotheses (L), (H0), (Ha), (Hab) and (Hb) hold. The \textit{BVH} hypothesis for the regrouped data $T_{j,n}'s$ is satisfied. The \textit{UAN} for regrouped data is ensured by the uniform negligibility of the variances (See  \eqref{BN-R}) which itself is forced by a condition on the uniform negligibility of the variances for non-regrouped data (See \eqref{BN-S}).\\

\Ni So hypotheses (L), (H0), (Ha), (Hab) and (Hb) hold and if $B^{\prime}_n \rightarrow 0$, the \textit{UAN} condition and the \textit{BVH} are satisfied and hence, the \textit{FLL} theorem for independent data $T_{j,n}$ applies and Lynderberg condition for regrouped data is forced by the Lynderberg-type condition for the non-regrouped data \eqref{LYND-NRD} (See \eqref{Lynder-O-2} and \eqref{Lynder-O-2}). Applying the Feller-L\'evy-Lynderberg (as Theorem 23 in \cite{ips-mfpt-ang}, Chapter 7, section 4) to get Part (2) of Theorem  \ref{general_WL_MPS_ASSOC}

\Ni Also for data having $2+\delta$ moments, for $\delta>0$, the Lyapounov condition for regrouped data is forced by the same condition for non-regrouped data \ref{L-NRD} (See \eqref{Lyap-RS} and \eqref{Lyap-S}). Hence the Lyapounov theorem for the $T_{j,n}$'s applies under the Lyapounov-type condition ($A^{\prime}_n \rightarrow 0$). We apply the Lyapounov theorem (as in in \cite{ips-mfpt-ang}) to get Part (1) of Theorem 22 \ref{general_WL_MPS_ASSOC}, Chapter 7, Section 4). $\blacksquare$\\

\section{Conclusion} \label{sec_05} 

\Ni The paper will be closed by Section \ref{sec_06} where the full computations of weak limits to the Gaussian law for moving partial sums of independent data are given. The paper, after offering a general handling moving partial sums convergence for independent data opens the rich field of such kinds of asymptotic weak laws for dependent data and their application to invariance principles and beyond to modeling in Applied statistics. The next step is to remain of these two kinds of dependence but to characterize asymptotic of random sums with not necessarily the Poissonian hypothesis on the counting process $N(\circ)$ in classical hypotheses.

\newpage
\section{Appendix} \label{sec_06}

\noindent \textbf{Proof of Theorem \ref{th_la}}. In the partial sums beginning by $p(n)=0$, by using Lemma 8 in \cite{ips-mfpt-ang},  Theorem \ref{th_la} holds for $\delta >1$ if it does for $\delta =1$. So it is enough to prove the theorem for $0<\delta \leq 1.$ \ By Lemma 9 in \cite{ips-mfpt-ang}, Assumption (\ref{cond_lya}) implies $s_{n}{\prime}\rightarrow +\infty $ and%
\begin{eqnarray*}
&&\max_{p(n)+1\leq k\leq n+p(n)}\left( \frac{\sigma _{k}}{s_{n}{\prime}}\right) ^{2+\delta }\leq
\max_{p(n)+1\leq k\leq n+p(n)}\frac{\mathbb{E}\left\vert X_{k}\right\vert ^{2+\delta }%
}{s_{n}^{\prime2+\delta }}\\
&&\leq \frac{1}{s_{n}^{\prime2+\delta }}\sum_{k=p(n)+1}^{n+p(n)}\mathbb{E}%
\left\vert X_{k}\right\vert ^{2+\delta }=A_{n}^{\prime}(\delta )\rightarrow 0.
\end{eqnarray*}

\bigskip \noindent Let us use the expansion of the characteristic functions 
\begin{equation*}
f_{k}(u)=\int e^{iux}dF_{k}(x)
\end{equation*}

\bigskip \noindent at the order two to get for each $k$, $p(n)+1\leq k\leq n+p(n)$, as
given in Lemma 6 in \cite{ips-mfpt-ang}, 
\begin{equation}
f_{k}(u/s_{n}{\prime})=1-\frac{u^{2}}{2}\frac{\sigma _{k}^{2}}{s_{n}^{\prime2}}+\theta
\frac{\left\vert u\right\vert ^{2+\delta }\mathbb{E}\left\vert
X_{k}\right\vert ^{2+\delta }}{s_{n}^{\prime2+\delta }},\, \, \left|\theta\right|<1.  \label{expan-eachk}
\end{equation}

\bigskip \noindent Now the characteristic function of $S_{n}{\prime}/s_{n}{\prime}$ is, for $%
u\in \mathbb{R},$
\begin{equation*}
f_{S_{n}{\prime}/s_{n}{\prime}}(u)=\prod\limits_{k=p(n)+1}^{n+p(n)}f_{k}(u/s_{n}{\prime}),
\end{equation*}

\noindent  that is
\begin{equation*}
\log f_{S_{n}{\prime}/s_{n}{\prime}}(u)=\sum_{k=p(n)+1}^{n+p(n)}\log f_{k}(u/s_{n}{\prime}).
\end{equation*}

\bigskip \noindent Now use the uniform expansion of $\log (1+u)$ at the
neighborhood at $0$, that is

\begin{equation}
\sup_{\left\vert u\right\vert \leq z}\left\vert \frac{\log (1+u)}{u}-1\right\vert =\varepsilon (z)\rightarrow 0.  \label{expans-log}
\end{equation}

\Bin  For each $k$ in (\ref{expan-eachk}), we have
\begin{equation}
f_{k}(u/s_{n}{\prime})=1-u_{kn}  \label{expans-unif}
\end{equation}

\bigskip \noindent with the uniform bound%
\begin{eqnarray*}
\left\vert u_{kn}\right\vert &\leq& \frac{\left\vert u\right\vert ^{2}}{2} \max_{p(n)+1\leq k\leq n+p(n)} 
\left(\frac{\sigma _{k}^{2}}{s_{n}^{\prime2}}\right)+\frac{\left\vert
u\right\vert ^{2+\delta }\mathbb{E}\left\vert X_{k}\right\vert ^{2+\delta }}{%
s_{n}^{\prime2+\delta }}\\
&=&\frac{\left\vert u\right\vert ^{2}}{2}\max_{p(n)+1\leq k\leq n+p(n)}\left( \frac{\sigma _{k}^{2}}{s_{n}^{\prime2}}\right) +\frac{\left\vert u\right\vert
^{2+\delta }\sum_{k=p(n)+1}^{n+p(n)}\mathbb{E}\left\vert X_{k}\right\vert ^{2+\delta }%
}{s_{n}^{\prime2+\delta }}=:u_{n}\rightarrow0.
\end{eqnarray*}

\bigskip \noindent Apply (\ref{expans-log}) to (\ref{expans-unif}) to have%
\begin{equation*}
\log f_{k}(u/s_{n}{\prime})=-u_{kn}+\theta _{n}u_{kn}\varepsilon (u_{n}), \, \, \left|\theta_{n}\right|<1
\end{equation*}

\bigskip \noindent and next 
\begin{eqnarray*}
\log f_{S_{n}{\prime}/s_{n}{\prime}}(u)&=&\sum_{k=p(n)+1}^{n+p(n)}\log f_{k}(u/s_{n}{\prime})\\
&=&-\frac{u^{2}}{2}+\left\vert u\right\vert ^{2+\delta }\theta A_{n}^{\prime}(\delta )+(\frac{u^{2}}{
2}-\left\vert u\right\vert ^{2+\delta }\theta A_{n}^{\prime}(\delta ))\varepsilon(u_{n})\theta _{n}\\
&\rightarrow& -u^{2}/2.
\end{eqnarray*}

\bigskip \noindent We get for u fixed,%
\begin{equation*}
f_{S_{n}{\prime}/s_{n}{\prime}}(u)\rightarrow \exp (-u^{2}/2).
\end{equation*}

\bigskip \noindent This completes the proof.\newline

\newpage
\bigskip  \noindent \textbf{Proof of Theorem \ref{th_LLF}}.\newline
\noindent The proof follows the lines of that of \cite{loeve}. But they are extended by more details and adapted and changed in some parts. Much
details were omitted. We get them back for making the proof understandable for students who just finished the measure and probability course.\newline

\bigskip \noindent Before we begin, let us establish an important property when (\ref{th_LLFb}) holds. Suppose that this latter holds.
We want to show that there exists a sequence $\varepsilon _{n}\rightarrow 0$
such that $\varepsilon _{n}^{-2}g_{n}(\varepsilon _{n})\rightarrow 0$ (this
implying $\ $also that $\varepsilon _{n}^{-1}g_{n}(\varepsilon
_{n})=o(\varepsilon _{n})\rightarrow 0$ and that $g_{n}(\varepsilon
_{n})=o(\varepsilon _{n}^{2})\rightarrow 0).$ To this end, let $k\geq 1$
fixed$.$ Since $g_{n}(1/k)\rightarrow 0$ as $n\rightarrow \infty ,$ we have $%
0\leq g_{n}(1/k)\leq k^{-3}$ for $n$ large enough.\newline

\bigskip \noindent We will get what we want from an induction on this
property. Fix $k=1$ and denote $n_{1}$ an integer such that $0\leq
g_{n}(1)\leq 1^{-3}$ for $n\geq n_{1}.$ Now we apply the same property on
the sequence $\{g_{n}(\circ ),n_{1}+1\}$ with $k=2$. We find a $n_{2}>n_{1}$
such that $0\leq g_{n}(1/2)\leq 2^{-3}$ for $n\geq n_{2}.$ Next we apply the
same property on the sequence $\{g_{n}(\circ ),n_{2}+1\}$ with $k=3$. We
find a $n_{3}>n_{2}$ such that $0\leq g_{n}(1/3)\leq 3^{-3}$ for $n\geq
n_{3}.$ Finally, there exists an infinite sequence of integers $%
n_{1}<n_{2}<...<n_{k}<n_{k+1}<...$ such that for each $k\geq 1,$ one has $%
0\leq g_{n}(1/k)\leq k^{-3}$ for $n\geq n_{k}.$

\bigskip \noindent Thus for each $n\geq 1$, for each $\forall k\geq 1$, there exists $k{(n)}$ such that 

$$
n_{k(n)}\leq n<n_{k(n)+1} \ and \ 0\leq g_{n}(1/k(n))\leq k(n)^{-3}.
$$

\Bin Put 

\begin{equation*}
\varepsilon _{n}=1/k(n) \text{ on }n_{k(n)}\leq n<n_{k(n)+1}
\end{equation*}

\bigskip \noindent We surely have $\varepsilon _{n}\rightarrow 0$ and $%
\varepsilon _{n}^{-2}g_{n}(\varepsilon _{n}).$ This is clear from%
\begin{equation*}
\left\{ 
\begin{tabular}{lll}
$\varepsilon _{n}=1/k(n)$ & on & $n_{k(n)}\leq n<n_{k(n)+1}$ \\ 
$\varepsilon _{n}^{-2}g_{n}(\varepsilon _{n})=k(n)^{2}(1/k(n)^{3})\leq (1/k(n))$ & on
& $n_{k(n)}\leq n<n_{k(n)+1}$%
\end{tabular}%
\right. .
\end{equation*}

\bigskip \noindent Now, we are going to use
\begin{equation}
\varepsilon _{n}\rightarrow 0\text{ and }\varepsilon
_{n}^{-2}g_{n}(\varepsilon _{n})\rightarrow 0.  \label{seqTravail}
\end{equation}

\bigskip \noindent \textbf{Proof\ of : (\ref{th_LLFb}) }$\Longrightarrow $%
\textbf{(\ref{th_LLFa}).}
\noindent Suppose (\ref{th_LLFb}) holds$.$ Thus there exists a
sequence $(\varepsilon _{n})_{n\geq 0}$ of positive numbers such that (\eqref{seqTravail}) prevails. First, we see that, for each $p(n)+1\leq k\leq n+p(n)$,

\begin{eqnarray*}
\frac{\sigma _{k}^{2}}{s_{n}^{\prime2}} &=&\frac{1}{s_{n}^{\prime2}}\int x^{2}dF_{k}(x)=%
\frac{1}{s_{n}^{\prime2}}\left\{ \int_{(\left\vert x\right\vert \geq \varepsilon
_{n}s_{n}{\prime})}x^{2}dF_{k}(x)+\int_{(\left\vert x\right\vert < \varepsilon
_{n}s_{n}{\prime})}x^{2}dF_{k}(x)\right\}  \\
&\leq &\frac{1}{s_{n}^{\prime2}}\int_{(\left\vert x\right\vert \geq \varepsilon
_{n}s_{n}{\prime})}x^{2}dF_{k}(x)+\varepsilon _{n}^{2} \\
&\leq &\frac{1}{s_{n}^{\prime2}}\sum_{k=p(n)+1}^{n+p(n)}\int_{(\left\vert x\right\vert \geq
\varepsilon _{n}s_{n}{\prime})}x^{2}dF_{k}(x)+\varepsilon_n^2=g(\varepsilon _{n})+\varepsilon _{n}^{2}.
\end{eqnarray*}

\bigskip \noindent It follows that
\begin{equation*}
\max_{p(n)+1\leq k\leq n+p(n)}\frac{\sigma _{k}^{2}}{s_{n}^{\prime2}}\leq g(\varepsilon
_{n})+\varepsilon _{n}^{2}\rightarrow \text{0.}
\end{equation*}

\bigskip \noindent Its remains to prove that $S_{n}{\prime}/s_{n}{\prime}$ $\rightsquigarrow
\mathcal{(0,1)}.$ To this end we are going to use this array of truncated randoms
variables $\{X_{nk}, \ p(n)+1\leq k\leq n+p(n), \ n\geq 1\}$ defined as follows. For each
fixed $n\geq 1,$ define

\begin{equation*}
X_{nk}=\left\{ 
\begin{tabular}{lll}
$X_{k}$ & if & $\left\vert X_{k}\right\vert \leq \varepsilon _{n}s_{n}^{\prime}$ \\ 
$0$ & if & $\left\vert X_{k}\right\vert >\varepsilon _{n}s_{n}^{\prime}$%
\end{tabular}%
\right. ,p(n)+1\leq k\leq n+p(n).
\end{equation*}

\bigskip \noindent Now consider summands $S_{nn}^{\prime}$ and $s_{nn}^{\prime2}$ defined by

$$
S_{nn}^{\prime}=\sum_{k=p(n)+1}^{n+p(n)} X_{nk} \,\,\,\,\, and \,\,\,\,\, s_{nn}^{\prime2}=\sum_{k=p(n)+1}^{n+p(n)} \mathbb{E}\left(X_{nk}-\mathbb{E}(X_{nk})\right)^2.
$$

\bigskip \noindent We remark that for any $\eta >0,$

\begin{eqnarray*}
\mathbb{P}\left( \left\vert \frac{S_{nn}^{\prime}}{s_{n}^{\prime}}-\frac{S_{n}^{\prime}}{s_{n}{\prime}}\right\vert >\eta
\right) \leq \mathbb{P}\left( \frac{S_{nn}^{\prime}}{s_{n}^{\prime}}\neq \frac{S_{n}^{\prime}}{s_{n}^{\prime}}\right)
\end{eqnarray*}

\bigskip \noindent and remark that
\begin{eqnarray*}
\left( \frac{S_{nn}^{\prime}}{s_{n}^{\prime}}\neq \frac{S_{n}^{\prime}}{s_{n}^{\prime}}\right) &=&\left( \exists
(p(n)+1\leq k\leq n+p(n)),X_{nk}\neq X_{k}\right)\\
&=&\left( \exists (p(n)+1\leq k\leq n+p(n)),\left\vert X_{k}\right\vert >\varepsilon _{n}s_{n}^{\prime}\right)
=\bigcup\limits_{k=p(n)+1}^{n+p(n)}(\left\vert X_{k}\right\vert >\varepsilon _{n}s_{n}^{\prime}).
\end{eqnarray*}

\bigskip \noindent We get
\begin{eqnarray*}
\mathbb{P}\left( \left\vert \frac{S_{nn}^{\prime}}{s_{n}{\prime}}-\frac{S_{n}^{\prime}}{s_{n}^{\prime}}\right\vert >\eta
\right)  &\leq &\sum_{k=p(n)+1}^{n+p(n)}\mathbb{P}(\left\vert X_{k}\right\vert >\varepsilon
_{n}s_{n}{\prime}) \\
&\leq &\sum_{k=p(n)+1}^{n+p(n)}\int_{(\left\vert x\right\vert \geq \varepsilon
_{n}s_{n}^{\prime})}dF_{k}(x)=\sum_{k=p(n)+1}^{n+p(n)}\int_{(\left\vert x\right\vert \geq
\varepsilon _{n}s_{n})}\left\{ \frac{1}{x^{2}}\right\} x^{2}dF_{k}(x) \\
&\leq &\left\{ \frac{1}{(\varepsilon _{n}s_{n}^{\prime})^{2}}\right\}
\sum_{k=p(n)+1}^{n+p(n)}\int_{(\left\vert x\right\vert \geq \varepsilon
_{n}s_{n}^{\prime})}x^{2}dF_{k}(x) \\
&\leq &\varepsilon _{n}^{-2}g_{n}(\varepsilon _{n})\rightarrow 0.
\end{eqnarray*}

\bigskip \noindent Thus $S_{nn}^{\prime}/s_{n}{\prime}$ and $S_{n}^{\prime}/s_{n}^{\prime}$ are equivalent in
probability. This implies that they have the same limit law or do not have a
limit law together. So to prove that $S_{n}^{\prime}/s_{n}^{\prime}$ has a limit law, we may
prove that $S_{nn}^{\prime}/s_{n}^{\prime}$ has a limit law. Next by Slutsky lemma, it will
suffice to establish the limit law of $S_{nn}^{\prime}/s_{nn}^{\prime}$
whenever we show that $s_{nn}^{\prime}/s_{n}^{\prime}\rightarrow 1.$ We focus on this. Begin
to remark that, since $\mathbb{E}(X_{k})=0$, we have the decomposition 
\begin{equation*}
0=\mathbb{E}(X_{k})=\int xdF_{k}(x)=\int_{(\left\vert x\right\vert \leq \varepsilon
_{n}s_{n}{\prime})}xdF_{k}(x)+\int_{(\left\vert x\right\vert >\varepsilon
_{n}s_{n}{\prime})}xdF_{k}(x)
\end{equation*}

\bigskip \noindent to get that 
\begin{equation*}
\left\vert \int_{(\left\vert x\right\vert \leq \varepsilon
_{n}s_{n}^{\prime})}xdF_{k}(x)\right\vert =\left\vert \int_{(\left\vert x\right\vert
>\varepsilon _{n}s_{n}^{\prime})}xdF_{k}(x)\right\vert .
\end{equation*}

\bigskip \noindent Remarking also that%
\begin{eqnarray}
\mathbb{E}(X_{nk})&=&\int_{(\left\vert X_{k}\right\vert \leq \varepsilon
_{n}s_{n}^{\prime})}X_{nk}d\mathbb{P}+\int_{(\left\vert X_{k}\right\vert >\varepsilon
_{n}s_{n}^{\prime})}X_{nk}d\mathbb{P} \label{trunc00} \\
&=&\int_{(\left\vert X_{k}\right\vert \leq \varepsilon
_{n}s_{n}^{\prime})}X_{k}dF_k+\int_{(\left\vert X_{k}\right\vert >\varepsilon
_{n}s_{n}^{\prime})}0dF_k.  \notag
\end{eqnarray}

\bigskip \noindent Combining all that leads to
\begin{eqnarray*}
\left\vert \mathbb{E}(X_{nk})\right\vert &=&\left\vert \int_{(\left\vert
X_{k}\right\vert \leq \varepsilon _{n}s_{n}^{\prime})}X_{k}d\mathbb{F_k}\right\vert\\
&=&\left\vert\int_{(\left\vert x\right\vert >\varepsilon _{n}s_{n}^{\prime})}xdF_{k}(x)\right\vert
\leq \int_{(\left\vert x\right\vert >\varepsilon _{n}s_{n}^{\prime})} \left|x\right|dF_{k}(x) \\
&= & \int_{(\left\vert x\right\vert >\varepsilon _{n}s_{n}^{\prime})}\frac{%
1}{\left\vert x\right\vert }x^{2}dF_{k}(x) \leq \frac{1}{\varepsilon
_{n}s_{n}^{\prime}} \int_{(\left\vert x\right\vert \geq\varepsilon
_{n}s_{n}^{\prime})}x^{2}dF_{k}(x).
\end{eqnarray*}

\bigskip \noindent Therefore%
\begin{equation}
\frac{1}{s_{n}^{\prime}}\sum_{k=p(n)+1}^{n+p(n)}\left\vert \mathbb{E}(X_{nk})\right\vert \leq
\varepsilon _{n}^{-1}g(\varepsilon _{n})\rightarrow 0.  \label{expect}
\end{equation}

\bigskip \noindent With help of this let us evaluate $s_{nn}^{\prime}/s_{n}^{\prime}.$ Notice
that for each fixed $n\geq 1,$ $\ $the $X_{nk}$ are still independent. The
technic used in (\ref{trunc00}) may be summarized as follows : any 
measurable function $g(\circ )$ such that $g(0)=0,$

\begin{equation} \label{ext_mes}
\mathbb{E}g(X_{nk})=\int_{(\left\vert X_{k}\right\vert \leq \varepsilon
_{n}s_{n}^{\prime})}g(X_{nk})d\mathbb{P}+\int_{(\left\vert X_{k}\right\vert >\varepsilon
_{n}s_{n}^{\prime})}g(0)d\mathbb{F_k}=\int_{(\left\vert X_{k}\right\vert \leq \varepsilon
_{n}s_{n}^{\prime})}g(X_{k})d\mathbb{F_k}.
\end{equation}

\bigskip \noindent By putting these remarks together, we obtain

\begin{eqnarray*}
1-\frac{s_{nn}^{\prime2}}{s_{n}^{\prime2}} &=&\frac{s_{n}^{\prime2}-s_{nn}^{\prime2}}{s_{n}^{\prime2}} \\
&=&\frac{1}{s_{n}^{\prime2}}\left\{
\sum\limits_{k=p(n)+1}^{n+p(n)}\mathbb{E}X_{k}^{2}-\sum\limits_{k=p(n)+1}^{n+p(n)}\mathbb{E}(X_{nk}-\mathbb{E}(X_{nk}))^{2}%
\right\}  \\
&=&\frac{1}{s_n^{\prime2}}\left\{ \sum\limits_{k=p(n)+1}^{n+p(n)}\mathbb{E}X_{k}^{2}-\left(
\sum\limits_{k=p(n)+1}^{n+p(n)}\mathbb{E}(X_{nk}^{2})-\mathbb{E}(X_{nk})^{2}\right) \right\}  \\
&=&\frac{1}{s_{n}^{\prime2}}\left\{
\sum\limits_{k=p(n)+1}^{n+p(n)}\mathbb{E}X_{k}^{2}-\sum\limits_{k=p(n)+1}^{n+p(n)}\mathbb{E}X_{nk}^{2}+\sum%
\limits_{k=p(n)+1}^{n+p(n)}\left( \mathbb{E}X_{nk}\right) ^{2}\right\}  \\
&=&\frac{1}{s_{n}^{\prime2}}\left\{ \sum\limits_{k=p(n)+1}^{n+p(n)}\int
X_{k}^{2}dF_k-\sum\limits_{k=p(n)+1}^{n+p(n)}\int_{(\left\vert X_{k}\right\vert \leq
\varepsilon _{n}s_{n}^{\prime})}X_{k}^{2}dF_k+\sum\limits_{k=p(n)+1}^{n+p(n)}\left( \mathbb{E}X_{nk}\right)
^{2}\right\}.  \\
\end{eqnarray*}

\noindent This leads to 

\begin{eqnarray*}
\vert1-\frac{s_{nn}^{\prime2}}{s_{n}^{\prime2}} \vert
&&=\frac{1}{s_{n}^{\prime2}}\left\{ \sum\limits_{k=p(n)+1}^{n+p(n)}\int_{(\left\vert X_{k}\right\vert >\varepsilon_{n}s_{n}^{\prime})}X_{k}^{2}dF_k+\sum\limits_{k=p(n)+1}^{n+p(n)}\left( \mathbb{E}X_{nk}\right)
^{2}\right\}  \\
&&\leq \frac{1}{s_{n}^{\prime2}}\left\{ \sum\limits_{k=p(n)+1}^{n+p(n)}\int_{(\left\vert
X_{k}\right\vert >\varepsilon
_{n}s_{n}^{\prime})}X_{k}^{2}dF_k+\sum\limits_{k=p(n)+1}^{n+p(n)}\left\vert \mathbb{E}
X_{nk} \right\vert ^{2}\right\} 
\end{eqnarray*}

\bigskip \noindent Finally, use the simple inequality of real numbers $%
\left( \sum_{i} \left\vert a_{i}\right\vert \right) ^{2}=\sum_{i} \left\vert
a_{i}\right\vert ^{2}+\sum_{i\neq j}\left\vert a_{i}\right\vert
\left\vert a_{j}\right\vert \geq \sum_{i} \left\vert a_{i}\right\vert ^{2}$ and
conclude from the last inequality that

\begin{eqnarray*}
\left\vert 1-\frac{s_{nn}^{\prime2}}{s_{n}^{\prime2}}\right\vert  &\leq &\frac{1}{s_{n}^{\prime2}}\left\{ \sum\limits_{k=p(n)+1}^{n+p(n)}\int_{(\left\vert X_{k}\right\vert
>\varepsilon _{n}s_{n}^{\prime})}X_{k}^{2}dF_k+\sum\limits_{k=p(n)+1}^{n+p(n)} \vert \mathbb{E}\left(X_{nk}\right)\vert^{2}\right\}  \\
&\leq &\frac{1}{s_{n^{\prime2}}}\left\{ \sum\limits_{k=p(n)+1}^{n+p(n)}\int_{(\left\vert
X_{k}\right\vert \geq\varepsilon _{n}s_{n}^{\prime})}X_{k}^{2}dF_k+\left(\sum\limits_{k=p(n)+1}^{n+p(n)}\vert\mathbb{E} X_{nk}\vert \right) ^{2}\right\}\\
&=&g(\varepsilon _{n})+\left( \frac{1}{s_{n}^{\prime}}\sum\limits_{k=p(n)+1}^{n+p(n)} \vert\mathbb{E} X_{nk}\vert \right)^2 .
\end{eqnarray*}

\bigskip \noindent By (\ref{expect}) above, we arrive at%
\begin{equation*}
\left\vert 1-\frac{s_{nn}^{\prime2}}{s_{n}^{\prime2}}\right\vert \leq g(\varepsilon
_{n})+\left(\varepsilon _{n}^{-1}g(\varepsilon _{n})\right)^2\rightarrow 0.
\end{equation*}

\bigskip \noindent It comes that $s_{nn}^{\prime}/s_{n}^{\prime}\rightarrow 1.$ Finally, the
proof of this part will derive  from the limit law of $S_{nn}^{\prime}/s_{nn}^{\prime}.$ We center the $X_{nk}$
at their expectations. To prove that the new summands $T_{nn}^{\prime}/s_{nn}^{\prime}$, where 

$$
T^{\prime}_{nn} = \sum_{k=p(n)+1}^{n+p(n)}\left(X_{nk} - \mathbb{E}\left(X_{nk}\right)\right),
$$ 

\Bin converge to $\mathcal{N}(0,1)$, we check the Lyapounov's condition for $\delta=1$, that is

\begin{equation*}
\frac{1}{s_{nn}^{\prime3}}\sum_{k=p(n)+1}^{n+p(n)}\mathbb{E}\left\vert
X_{nk}-\mathbb{E}X_{nk}\right\vert ^{3}\rightarrow 0\text{ as }n\rightarrow \infty .
\end{equation*}

\bigskip \noindent But we have
\begin{eqnarray*}
\frac{1}{s_{nn}^{\prime3}}\sum_{k=p(n)+1}^{n+p(n)}\mathbb{E}\left\vert X_{nk}-\mathbb{E}X_{nk}\right\vert ^{3} &=&\frac{1}{s_{nn}^{\prime3}}\sum_{k=p(n)+1}^{n+p(n)}\mathbb{E}\biggr(\left\vert X_{nk}-\mathbb{E}X_{nk}\right\vert \times \left\vert X_{nk}-\mathbb{E}X_{nk}\right\vert ^{2}\biggr) \\
&\leq &\frac{\mathbb{E}\vert X_{nk}\vert}{s_{nn}^{\prime3}}\sum_{k=p(n)+1}^{n+p(n)} \mathbb{E}\left\vert X_{nk}-\mathbb{E}X_{nk}\right\vert ^{2} \\
&+&\frac{2}{s_{nn}^{\prime3}}\sum_{k=p(n)+1}^{n+p(n)}\mathbb{E}\left(\left\vert X_{nk}\right\vert \left\vert X_{nk}-\mathbb{E}X_{nk}\right\vert ^{2}\right).
\end{eqnarray*}

\bigskip \noindent Take $g(\circ )=\left\vert \circ \right\vert $ in (\ref{ext_mes}) to see again that

\begin{equation*}
\mathbb{E}\left\vert X_{nk}\right\vert =\int_{(\left\vert X_{k}\right\vert
\leq \varepsilon _{n}s_{n}^{\prime})}\left\vert X_{k}\right\vert dF_k\leq \varepsilon
_{n}s_{n}^{\prime}.
\end{equation*}

\noindent We also have

\begin{eqnarray*}
\mathbb{E}\left(\left\vert X_{nk}\right\vert \left\vert X_{nk}-\mathbb{E}X_{nk}\right\vert ^{2}\right)
&=&\int_{(|X_k|\leq \epsilon_n s_{n}^{\prime})} \left\vert X_{nk}\right\vert \left\vert X_{nk}-\mathbb{E}X_{nk}\right\vert ^{2} \ d\mathbb{P}\\
&\leq & \varepsilon_{n}s_{n}^{\prime} \int_{(|X_k|\leq \epsilon_n s_{n}^{\prime})} \left\vert X_{nk}-\mathbb{E}X_{nk}\right\vert ^{2} \ d\mathbb{P}.
\end{eqnarray*}

\bigskip \noindent The last formula yield
\begin{eqnarray*}
\frac{1}{s_{nn}^{\prime3}}\sum_{k=p(n)+1}^{n+p(n)}\mathbb{E}\left\vert
X_{nk}-\mathbb{E}X_{nk}\right\vert ^{3}&\leq& \frac{2\varepsilon _{n}s_{n}^{\prime}}{s_{nn}^{\prime3}}%
\sum_{k=p(n)+1}^{n+p(n)}\mathbb{E}\left\vert X_{nk}-\mathbb{E}X_{nk}\right\vert ^{2}\\
&=&\frac{2\varepsilon _{n}s_{n}^{\prime}s_{nn}^{\prime2}}{s_{nn}^{\prime3}}=\varepsilon _{n}\frac{2s_{n}^{\prime}}{%
s_{nn}^{\prime}}\rightarrow 0.
\end{eqnarray*}

\bigskip \noindent By Lyapounov's theorem, we have
\begin{equation*}
\frac{T_{nn}^{\prime}}{s_{nn}^{\prime}}=\frac{S_{nn}^{\prime}-\sum_{k=p(n)+1}^{n+p(n)}\mathbb{E}(X_{nk})}{s_{nn}^{\prime}}%
\rightsquigarrow \mathcal{N}(0,1).
\end{equation*}

\bigskip \noindent Since $s_{nn}^{\prime}/s_{n}^{\prime}\rightarrow 1$ and by (\ref{expect})%
\begin{equation*}
\left\vert \frac{\sum_{k=p(n)+1}^{n+p(n)}\mathbb{E}(X_{nk})}{s_{nn}^{\prime}}\right\vert \leq \frac{s_{n}^{\prime}%
}{s_{nn}^{\prime}}\left\{ \frac{1}{s_{n}^{\prime}}\sum_{k=p(n)+1}^{n+p(n)}\vert\mathbb{E} X_{nk}\vert
)\right\} \rightarrow 0.
\end{equation*}

\bigskip \noindent We conclude that $S_{nn}^{\prime}/s_{nn}^{\prime}$ converges to $\mathcal{N}(0,1)$.%
\newline

\bigskip \noindent \textbf{Proof of : (\ref{th_LLFa})$\Longrightarrow$ (\ref{th_LLFb})}.
The convergence to $\mathcal{N}(0,1)$ implies that for any fixed $t\in \mathbb{R},$ we have%
\begin{equation}
\prod\limits_{k=p(n)+1}^{n+p(n)}f_{k}(u/s_{n}^{\prime})\rightarrow \exp (-u^{2}/2).
\label{limit}
\end{equation}

\bigskip \noindent We are going to use uniform expansions of $\log (1+z)$.
We have 
\begin{equation*}
\lim_{z\rightarrow 0}\left\vert \frac{\log (1+z)-z}{z^{2}}\right\vert =\frac{%
1}{2},
\end{equation*}

\bigskip \noindent this implies
\begin{equation}
\sup_{z\leq u}\left\vert \frac{\log (1+z)-z}{z^{2}}\right\vert =\varepsilon
(u)\rightarrow 1/2\text{ as }u\rightarrow 0.  \label{II_02}
\end{equation}

\bigskip \noindent Now, we use the expansion

\begin{equation*}
f_{k}(u/s_{n}^{\prime})=1+\theta _{k}\frac{u^{2}\sigma _{k}^{2}}{2s_{n}^{\prime2}}, \ \vert\theta_k\vert\leq 1
\end{equation*}

\bigskip \noindent to get that

\begin{equation}
\max_{p(n)+1\leq k\leq n+p(n)}\left\vert f_{k}(u/s_{n}^{\prime})-1\right\vert \leq \frac{%
u^{2}}{2}\max_{p(n)+1\leq k\leq n+p(n)}\frac{\sigma _{k}^{2}}{s_{n}^{\prime2}}%
=:u_{n}\rightarrow 0  \label{II_01}
\end{equation}

\noindent and next 

\begin{equation*}
\left\vert f_{k}(u/s_{n}^{\prime})-1\right\vert^2 =\theta _{k}^{2}\frac{u^{2}\sigma
_{k}^{2}}{2s_{n}^{\prime2}}\times \frac{u^{2}\sigma _{k}^{2}}{2s_{n}^{\prime2}}\leq %
\left[ \frac{u^{4}}{4}\max_{p(n)+1\leq k\leq n+p(n)}\frac{\sigma _{k}^{2}}{s_{n}^{\prime2}}%
\right] \times \frac{\sigma _{k}^{2}}{s_{n}^{\prime2}}.
\end{equation*}

\bigskip \noindent This latter implies

\bigskip 
\begin{equation*}
\sum\limits_{k=p(n)+1}^{n+p(n)}\left\vert f_{k}(u/s_{n}^{\prime})-1\right\vert^2 \leq \left[ 
\frac{u^{4}}{4}\max_{p(n)+1\leq k\leq n+p(n)}\frac{\sigma _{k}^{2}}{s_{n}^{\prime2}}\right]
=:B_{n}(u)\rightarrow 0.
\end{equation*}

\bigskip \noindent By (\ref{II_01}), we see that $\log f_{k}(u/s_{n}^{\prime}$%
\bigskip ) uniformly defined in $k$ such that $p(n)+1\leq k\leq n+p(n)$ for $n$ large enough and (%
\ref{limit}) becomes
\begin{equation*}
\sum\limits_{k=p(n)+1}^{n+p(n)}\log f_{k}(u/s_{n}^{\prime})\rightarrow -u^{2}/2
\end{equation*}

\bigskip \noindent that is 
\begin{equation*}
\frac{u^{2}}{2}+\sum\limits_{k=p(n)+1}^{n+p(n)}\log f_{k}(u/s_{n}^{\prime})\rightarrow 0.
\end{equation*}

\bigskip \noindent Now using the uniform bound of $\left\vert
f_{k}(u/s_{n}^{\prime})-1\right\vert $ by $u_{n}$ to get%
\begin{equation*}
\log f_{k}(u/s_{n}^{\prime})=f_{k}(u/s_{n}^{\prime})-1+\theta_{kn}(f_{k}(u/s_{n}^{\prime})-1)^{2}\varepsilon
(u_{n}), \ \left|\theta_{kn}\right|\leq 1
\end{equation*}

\bigskip \noindent and then
\begin{eqnarray*}
\frac{u^{2}}{2}+\sum\limits_{k=p(n)+1}^{n+p(n)}\log f_{k}\left(u/s_{n}^{\prime}\right) &=&\frac{u^{2}}{2}%
+\sum\limits_{k=p(n)+1}^{n+p(n)}f_{k}(u/s_{n}^{\prime})-1+\theta_{kn}(f_{k}(u/s_{n}^{\prime})-1)^{2}\varepsilon
(u_{n}) \\
&=&\left\{ \frac{u^{2}}{2}-\sum\limits_{k=p(n)+1}^{n+p(n)}\left(1-f_{k}(u/s_{n}^{\prime})\right)\right\}\\
&+&\left\{ \sum\limits_{k=p(n)+1}^{n+p(n)}\theta_{kn}(f_{k}(u/s_{n}^{\prime})-1)^{2}\right\} \varepsilon (u_{n}),
\end{eqnarray*}

\bigskip \noindent with
\begin{equation*}
\left\vert \left\{ \sum\limits_{k=p(n)+1}^{n+p(n)}\theta_{kn}(f_{k}(u/s_{n}^{\prime})-1)^{2}\right\}
\varepsilon (u_{n})\right\vert \leq B_{n}(u)\left\vert \varepsilon
(u_{n})\right\vert =o(1).
\end{equation*}

\bigskip \noindent We arrive at%
\begin{equation*}
\frac{u^{2}}{2}=\sum\limits_{k=p(n)+1}^{n+p(n)}\left(1-f_{k}(u/s_{n}^{\prime})\right)+o(1).
\end{equation*}

\bigskip \noindent If we take the real parts, we have for any fixed $\varepsilon >0$,

\begin{eqnarray*}
\frac{u^{2}}{2} &=&\sum\limits_{k=p(n)+1}^{n+p(n)}\int \left(1-\cos \frac{ux}{s_{n}^{\prime}}%
\right)dF_{k}(x)+o(1) \\
&=&\sum\limits_{k=p(n)+1}^{n+p(n)}\int_{(\left\vert x\right\vert <\varepsilon
s_{n}^{prime})}\left(1-\cos \frac{ux}{s_{n}^{\prime}}\right)dF_{k}(x)\\
&+&\sum\limits_{k=p(n)+1}^{n+p(n)}\int_{(\left\vert x\right\vert \geq \varepsilon s_{n}^{\prime})}\left(1-\cos \frac{ux}{s_{n}^{\prime}}\right)dF_{k}(x)+o(1),
\end{eqnarray*}

\bigskip \noindent that is
\begin{equation*}
\frac{u^{2}}{2}-\sum\limits_{k=p(n)+1}^{n+p(n)}\int_{(\vert x\vert <\varepsilon s_{n}^{\prime})}\left(1-\cos \frac{ux}{s_{n}^{\prime}}\right)dF_{k}(x)=\sum\limits_{k=p(n)+1}^{n+p(n)} \int_{(\vert x\vert \geq \varepsilon s_{n}^{\prime})}\left(1-\cos \frac{ux}{s_{n}^{\prime}}%
\right)dF_{k}(x)+o(1).
\end{equation*}

\bigskip \noindent We are going to use the following fact: for all $\delta ,0<\delta \leq 1$,

$$\sqrt{2(1-\cos a)}\leq2\left\vert a/2\right\vert ^{\delta }.$$

\Bin Apply this for $\delta =1$ to have%
\begin{eqnarray*}
&&\sum\limits_{k=p(n)+1}^{n+p(n)}\int_{(\vert x\vert <\varepsilon
s_{n}^{\prime})}\left(1-\cos \frac{ux}{s_{n}^{\prime}}\right)dF_{k}(x)\\
&\leq &\frac{u^{2}}{2s_{n}^{\prime2}}\sum\limits_{k=p(n)+1}^{n+p(n)}\int_{(\left\vert x\right\vert <\varepsilon
s_{n}^{\prime})}x^{2}dF_{k}(x) \\
&=&\frac{u^{2}}{2s_{n}^{\prime2}} \left\{\sum\limits_{k=p(n)+1}^{n+p(n)}\int x^{2}dF_{k}(x)-\sum\limits_{k=p(n)+1}^{n+p(n)}\int_{(\left\vert x\right\vert \geq
\varepsilon s_{n}^{\prime})}x^{2}dF_{k}(x)\right\}  \\
&=&\frac{u^{2}}{2s_{n}^{\prime2}} \left\{s_{n}^{\prime2}-\sum\limits_{k=p(n)+1}^{n+p(n)}\int_{(\left\vert
x\right\vert \geq \varepsilon s_{n}^{\prime})}x^{2}dF_{k}(x)\right\}=\frac{u^{2}}{2}%
(1-g_{n}(\varepsilon )).
\end{eqnarray*}

\bigskip \noindent On the other hand%
\begin{eqnarray*}
\sum\limits_{k=p(n)+1}^{n+p(n)}\int_{(\left\vert x\right\vert \geq \varepsilon
s_{n}^{\prime})}\left(1-\cos \frac{ux}{s_{n}^{\prime}}\right)dF_{k}(x) &\leq
&2\sum\limits_{k=p(n)+1}^{n+p(n)}\int_{(\left\vert x\right\vert \geq \varepsilon
s_{n}^{\prime})}dF_{k}(x) \\
&=&2\sum\limits_{k=p(n)+1}^{n+p(n)}\int_{(\left\vert x\right\vert \geq \varepsilon
s_{n}^{\prime})}\left\{ \frac{1}{x^{2}}\right\} x^{2}dF_{k}(x) \\
&\leq &\frac{2}{\varepsilon _{{}}^{2}s_{n}^{\prime2}}\sum\limits_{k=p(n)+1}^{n+p(n)}\int_{%
(\left\vert x\right\vert \geq \varepsilon s_{n}^{\prime})}x^{2}dF_{k}(x)\leq \frac{2}{%
\varepsilon ^{2}}.
\end{eqnarray*}

\bigskip \noindent By putting all this together, we have

\begin{equation*}
\frac{u^{2}}{2}\leq \frac{u^{2}}{2}(1-g_{n}(\varepsilon ))+\frac{2}{%
\varepsilon ^{2}}+o(1)
\end{equation*}

\bigskip \noindent which leads 

\begin{equation*}
\frac{u^{2}}{2}g_{n}(\varepsilon)\leq \frac{2}{\varepsilon ^{2}}+o(1)
\end{equation*}

\bigskip \noindent which in turns implies 
\begin{equation*}
g_{n}(\varepsilon )\leq \frac{2}{u^{2}}\left(\frac{2}{\varepsilon ^{2}}+o(1)\right) \ \ for \ all \ u\in \mathbb{R}^{\ast} \ and \ all \ n\geq 1.
\end{equation*}

\Bin So 
 
\begin{equation*}
\limsup_{n\rightarrow +\infty}\ g_{n}(\varepsilon )\leq \frac{2}{u^{2}}\left(\frac{2}{\varepsilon ^{2}}+o(1)\right) \ \ for \ all \ u\in \mathbb{R}^{\ast}. 
\end{equation*}

\bigskip \noindent By letting  $u\rightarrow +\infty$, we get 
\begin{equation*}
g_{n}(\varepsilon )\rightarrow 0.
\end{equation*}

\bigskip \noindent This concludes the proof. $\blacksquare$
\end{document}